\documentclass[12pt]{article}
\usepackage{amsmath, amssymb}
\usepackage{amsfonts}
\usepackage[dvips]{graphicx}

\newtheorem{theorem}{Theorem}[section]
\newtheorem{lemma}{Lemma}[section]
\newtheorem{definition}{Definition}[section]
\newtheorem{corollary}{Corollary}[section]
\newtheorem{remark}{Remark}[section]

\newtheorem{proposition}{Proposition}[section]

\def\QED{\mbox{\rule[0pt]{1.5ex}{1.5ex}}}

\def\endproof{\hspace*{\fill}~\QED\par\endtrivlist\unskip}
\def\keywords{\vspace{-.3em}
    \if@twocolumn
      \small\it Keywords\/\bf---$\!$%
    \else
      \begin{center}\small\bf Keywords\end{center}\quotation\small
    \fi}
\def\endkeywords{\vspace{0.6em}\par\if@twocolumn\else\endquotation\fi
    \normalsize\rm}

\def\Label#1{\label{#1}\ [\ #1 \ ]\ }
\def\Label{\label}

\title{General theory for integer-type algorithm \\ for higher order differential equations}

\author{Fuminori SAKAGUCHI\thanks{Faculty of Engineering, University of Fukui, 3-9-1 Bunkyo, Fukui 910-8507, Japan ({\tt fsaka@u-fukui.ac.jp})} \and Masahito HAYASHI\thanks{Graduate School of Information Sciences, Tohoku University, Sendai 980-8579, Japan ({\tt hayashi@math.is.tohoku.ac.jp}) Centre for Quantum Technologies, National University of Singapore, 3 Science Drive 2, Singapore 117542}}

\date{}
\begin{document}

\maketitle

\begin{abstract}
Based on functional analysis, we propose an algorithm 
for finite-norm solutions of higher-order linear 
Fuchsian-type 
ordinary differential equations (ODEs) $P(x,\frac{d}{dx})f(x)=0$ with 
$P(x,\frac{d}{dx}):=\displaystyle\Bigl( \sum_{m=0}^M p_m (x) \bigl({\textstyle \frac{d}{dx}}\bigr)^m \Bigr)$
by using only the four arithmetical operations on integers. 

This algorithm is based on a band-diagonal matrix representation of 
the differential operator $P(x,\frac{d}{dx})$,  
though it is quite different from the usual Galerkin methods.
This representation is made for 
the respective CONSs of the input Hilbert space ${\cal H}$ and the output Hilbert space 
${{\cal H}^\Diamond}$ of  $P(x,\frac{d}{dx})$.
This band-diagonal matrix enables the construction of a recursive algorithm for solving the ODE. 
However, a solution of the simultaneous linear equations represented by this matrix 
does not necessarily 
correspond to the true solution of ODE.
We show that
when this solution is an $\ell ^2$ sequence, it corresponds to the true solution of ODE.
We invent a method based on an integer-type algorithm for extracting only $\ell ^2$ components.
Further, the concrete choice of Hilbert spaces ${\cal H}$ and 
${{\cal H}^\Diamond}$ is also given for our algorithm
when $p_m$ is a polynomial or a rational function with rational coefficients.
We check how our algorithm works based on several numerical demonstrations 
related to special functions, where the results show that the accuracy of our method is
extremely high. 
\end{abstract}


AMS: 65L99, 42C15, 65L60, 34A45

\pagestyle{myheadings}
\thispagestyle{plain}

\section{Introduction}
\Label{sec:in}
Linear ordinary differential equations (ODE) of the type 
\begin{eqnarray}
P(x,{\textstyle \frac{d}{dx}})f(x):= \left( \sum_{m=0}^M p_m (x) \left( \frac{d}{dx}\right)^m \right) f(x) = 0 \Label{eqn:Haya-1}
\end{eqnarray}
are very important tools in many fields (physics, engineering etc.). 
In many useful cases, the functions $p_m (x)$ are polynomials or rational functions. 
As is well known, it is difficult in general to solve them analytically for higher-order cases 
although there are relatively general methods for second-order equations with low-degree polynomials $p_m (x)$, for which we employ hypergeometric functions (or special functions)~\cite{Eld} and power series expansions about nonsingular points or 
regular singular points~\cite{CoLe}.
(Practically, instead of analytical methods, many kinds of numerical methods have been proposed and used.)  
The aim of this paper is to obtain solutions $f(x)$ of a linear ODE (\ref{eqn:Haya-1}) in a Hilbert space ${\cal H}$ of functions on $\mathbb{R}$
when the equation is of higher order and/or the function $p_m (x)$ is of higher degree. 
Solutions with finite norm are sometimes very important in quantum mechanics (e.g. wavefunctions of particles bound by potentials)~\cite{Mes}, 
and for transit or temporary phenomena 
in signal processing and circuit theory 
that are almost localized in the time coordinate in many applications, for example. 

In this paper, we propose an integer-type general algorithm for solving these ODEs, 
by choosing function spaces and their basis systems appropriately. 
This method is based on a pair of Hilbert spaces ${\cal H}$ and ${{\cal H}^\Diamond}$ with distinct inner products, where the domain of the differential operator $P(x,{\textstyle \frac{d}{dx}})$ is a dense subspace of ${\cal H}$ and its range is a subspace of ${{\cal H}^\Diamond }$. 
Under appropriate choice of these spaces and their basis systems which will be presented in this paper, a differential equation can be expressed by band-diagonal-type simultaneous linear equations, and all the `matrix elements' are rational-(complex-)valued. 
Moreover, under the same choice, all the basis functions are rational functions.
In addition, from the properties of the basis functions used, this method has a somewhat similar feature to power series expansions about nonsingular points 
or regular singular points, in the sense that the solution can be expanded as linear combinations of the powers of a rational function of $x$ with rational coefficients. 
From another point of view, this method is closely related to the Laurent expansion and hyperfunctions in complex analysis and to the Fourier series, under some changes of variable.
Therefore, this method has a `semi-analytical' character and it can be discussed from the standpoint of mathematical analysis, though it is a kind of numerical method. 

Since it is difficult to apply analytical methods to general higher-order linear differential equations, 
various kinds of numerical methods have been proposed. One group is based on the discretization 
of coordinate or on the  
differences or on the relations between adjacent lattice points (Runge-Kutta methods, for example). Another group is based on finite-dimensional subspaces of an infinite-dimensional function space, such as 
the collocation method, the Ritz-Galerkin method and the Petrov-Galerkin method~\cite{KrVa}~\cite{Bre}, for example. 
In this group, many kinds of finite element methods~\cite{Bre}~\cite{ErGu} have been proposed and used widely and efficiently in many fields. 
These methods construct subspaces spanned by finite elements with very localized compact supports. 
In addition, this group contains a subgroup which uses subspaces spanned by {\it globally smooth} basis functions~\cite{KrVa}  
such as the Hermite functions. 

The method to be proposed in this paper 
is similar to the latter subgroup in the sense that it is based on the finite-dimensional subspaces spanned by global smooth basis functions. 

The proposed method is different from Ritz-Galerkin method, in that 
the function space ${{\cal H}^\Diamond}$ is different from ${\cal H}$, 
and ${{\cal H}^\Diamond}$ is wider than ${\cal H}$ in the proposed method. 
Here, remember that the function space ${\cal H}$ corresponds to the domain of 
the differential operator $P(x,{\textstyle \frac{d}{dx}})$
and the function space ${{\cal H}^\Diamond}$ does to the range of $P(x,{\textstyle 
\frac{d}{dx}})$.
The choice of different function spaces ${\cal H}$ and ${{\cal H}^\Diamond}$ may be possible even for Petrov-Galerkin method. 
However, the method proposed here is quite different from the usual `standard truncation methods' or `projection methods' such as the Ritz-Galerkin and Petrov-Galerkin methods, in respect of the following point: 


Although the 
Ritz-Galerkin and Petrov-Galerkin 
methods are based on the solutions of simultaneous linear equations with a square matrix truncated within a finite dimension,
the method proposed here 
is based on finite-dimensional truncations of the exact solutions of the {\it infinite-dimensional} simultaneous linear equations.  
There is a possibility that our numerical solution coincides exactly with the orthogonal projection of the true solution to the finite-dimensional subspace, whatever its dimension may be, 
and we have already had some numerical examples where this perfect coincidence occurs really~\cite{paper3}. 
In order to realize this direction, we solve 
simultaneous linear equations with a non-square-type band-diagonal matrix, in which its column is larger than its row.
Since the solution is not unique because this non-square-type band-diagonal matrix has a non-trivial kernel, 
we have to extract one solution among the above solutions.
In order to resolve this problem, we extract one solution among them using a novel method, which will be explained in the latter part of this introduction. Further, this matrix elements do not change when the dimension of the subspace increases. 
Hence the proposed method provides a recursive algorithm with no round-off errors up to an arbitrary dimension
for the vectors of the space of exact solutions of the infinite-dimensional simultaneous linear equations.

The method to be proposed has five advantages other than the integer-type property mentioned above, as follows; 
(1) especially when 
the ODE has no singular point 
or when the ODE belongs to the Fuchsian class even if 
it has singular points,  
this method can determine the structure itself of the function space of solutions in ${\cal H}$ of the differential equation, directly from the numerical results. 
(2) Another advantage is that the convergence of the error to $0$ is guaranteed as the dimension of the subspace tends to infinity, and an upper bound of the error can be given for the finite-dimensional case. 
(3) Moreover, it does not require any calculation of large  
matrices (inverse matrix, eigenvector, etc.) for solving our simultaneous linear equations. 
(4) Another strong point is that the basis functions of ${\cal H}$ are smooth sinusoidal-like wavepackets with spindle-shaped envelopes, which are suitable for the expansion
of various kind of `natural' functions decaying as $x\to\pm\infty $. 
In this sense, the basis functions contain both global and local information.
(5) Another strong point is that this method requires a small amount of calculations for obtaining high-accuracy solutions. For example, when the coefficients in the expansion of a true solution by the basis functions decay exponentially, the amount of calculations required by this method is almost proportional asymptotically to the cube of the number of required significant digits.

In this paper, we will show the validity of the band-diagonal matrix representation, i.e., 
we will show that the square-summable solutions of the simultaneous linear equations according to the band-diagonal matrix always correspond to true solutions of the ODE except at 
the singular points of the ODE. 
Especially, when 
the ODE has no singular point, 
we will show 
the one-to-one correspondence between the true solutions in ${\cal H}$ of the differential equation and the square-summable number sequences satisfying the simultaneous linear equations represented by the band-diagonal matrix. 
The larger part of its proofs is based mainly on functional analysis. 
Similar one-to-one correspondence can be proved for the cases of the Fuchsian class by a  modification even if 
the ODE has singular points. 

However, the presented method has 
a pitfall based on finite-dimensional truncations of the exact solutions of the infinite-dimensional simultaneous linear equations. 
This pitfall is due to the non-uniqueness of non-square-summable solutions of the simultaneous linear equations
because the number of linearly independent solutions of the simultaneous linear equation is not smaller than the bandwidth whereas the number of linearly independent solutions of the differential equations is not greater than its order $M$. 
That is, there are solutions which do not correspond to true solutions of the differential equations,
we call these solutions {\em extra} solutions.
In order to resolve this problem, we propose a method to to remove the extra solutions effectively. 
This method is based on quasi-minimization of the ratio between a norm sensitive to divergence and another norm insensitive to divergence. 
This quasi-minimization guarantees the convergence to $0$ of the error in the numerical results; 
the accuracy of the numerical results is sufficiently high, even for finite dimensions.

For minimization of the ratio between two quadratic forms,
the usual method is based on 
the eigenspace of the matrix $A^{-\frac{1}{2}}BA^{-\frac{1}{2}}$ with the two corresponding inner-product-matrices $A$ and $B$.
However, it is difficult to apply this method to the above problem, due to round-off errors, because these inner-product matrices are usually very close to a singular matrix with rank $1$. 
However, in this paper, we propose an alternative integer-type method for quasi-minimization, which does not require as much calculation as the usual method. 
This method is based on a kind of quasi-orthogonalization of integer-valued vectors, which is realized by an idea that is conceptually between the Gram-Schmidt process and the Euclidean algorithm. 

An integer-type recursive algorithm similar to the proposed one may be applied also to the Petrov-Galerkin method, in order to calculate the head and intermediate rows of solution vectors of the system of simultaneous linear equations described by a large-dimensional band-diagonal  square matrix. However, in this case, we have to calculate new linear combinations of solution vectors satisfying the final constraints (linear equations in the bottom rows). For this calculation, we can solve another system of simultaneous linear equations, which is described by another square matrix whose dimension is half a band width. However, this matrix is usually very close to a singular matrix of rank 1. Moreover, the elements in the bottom rows of the solution vectors are rational numbers whose numerators and denominators are huge integers. Therefore, this type of Petrov-Galerkin method requires a much larger amount of calculations than the proposed method based on non-square matrix and integer-type quasi-orthogonalization.

\begin{table}[bht]
\caption{Difference from Galerkin methods for norm-finite solutions (without initial conditions)}
\begin{center}
\begin{tabular}{|@{}c@{}|c|c|c@{}|}
\hline 
 Method & Ritz-Galerkin & Petrov-Galerkin& proposed \\
\hline 
 ${{\cal H}^\Diamond}$ and ${\cal H}$
& same & same/different & different \\
\hline 
Corresponding matrix & \multicolumn{2}{c|}{ square } & non-square and \\
with truncation 
& \multicolumn{2}{@{}c|}{(band-diagonal / general)} & band-diagonal \\
\hline 
Extra solutions
& \multicolumn{2}{@{}c|}{ No} & can be removed \\
\hline 
&
\multicolumn{2}{@{}c|}{ Eigenvalue and eigenvector} & 
Exact kernel vector  \\ 
Solution vector &
\multicolumn{2}{@{}c|}{ of finite-dimensional matrix }& (infinite-dimensional)
  \\ 
\hline
\end{tabular}
\end{center}
\label{tbl:diff}
\end{table}

The contents of the paper are as follows;
Section \ref{sec:ab} explains an abstract framework for a general algorithm,  
using a pair of Hilbert spaces ${\cal H}$ and ${{\cal H}^\Diamond }$ with distinct inner products. 
Subsection \ref{subsec:s21} states 
the basic conditions for the pair of Hilbert spaces ${\cal H}$ and ${{\cal H}^\Diamond }$.
This subsection gives a sufficient condition
for band-diagonal matrix representation for the ODE.
In Subsection \ref{subsec:s22}, 
using this band-diagonal form,
we provide the basic structure of our recursive algorithm 
with removal of non-$\ell^2$-components.
Section 3 presents the 
concrete choices of 
Hilbert spaces and basis systems used in our algorithm,
and checks that they satisfy the conditions given in Section \ref{sec:ab}.
Section 4 gives proofs of theorems mentioned in Section \ref{sec:ab}.
In Section 5, we give 
some numerical examples related to special functions 
and show how effectively our 
algorithm works, 
where we are successful to solve ODEs in a very high accuracy with a relatively small amount of calculations (approximately proportional to a power of the number of required significant digits, empirically). Moreover, we show how numerical results are successful even for the cases 
with singular points. 
Section 6 discusses a related topic and further extensions of our algorithm.

\section{Abstract structure of our algorithm}
\Label{sec:ab}

\begin{table}[bht]
\caption{Components of our algorithm}
\begin{center}
\begin{tabular}{|@{}c@{}|l@{}|@{}c@{}|}
\hline 
$P(x,{\textstyle \frac{d}{dx}})$ & Differential operator & Top of \ref{subsec:s21} \\
$M$ & Order of $P(x,{\textstyle \frac{d}{dx}})$ & Top of \ref{subsec:s21} \\
$p_m(x)$ & Polynomial: $m$-th order coefficient func. of $P(x,{\textstyle \frac{d}{dx}})$ & Top of \ref{subsec:s21} \\
${\cal H}$ & Input Hilbert space & Top of \ref{subsec:s21} \\
${{\cal H}^\Diamond }$ & Output Hilbert space & Top of \ref{subsec:s21} \\ 
$\tilde{A}_P$ & Operator ${\cal H}\to {\cal H}$ defined as action of $P(x,\frac{d}{dx})$ & Top of \ref{subsec:s21} \\
$A_P$ &Closed extension of $\tilde{A}_P$ & Top of \ref{subsec:s21}\\
$\tilde{B}_P$ & Operator ${\cal H}\to {{\cal H}^\Diamond }$ defined as action of $P(x,\frac{d}{dx})$ & Top of \ref{subsec:s21}\\
$B_P$ &Closed extension of $\tilde{B}_P$ & Top of \ref{subsec:s21}\\ 
$e_n$ & Basis of ${\cal H}$ &{\bf C1} \\ 
$e_n^\Diamond $ & Basis of ${\cal H}^\Diamond $ &{\bf C1} \\ 
$f_n $& $f_n:= \langle f,e_n \rangle _{{{\cal H}}}$ for $f\in {\cal H}$& Theorem \ref{thm:c:matrixeq}\\
$\vec{f}$ & Vector representation of $\{f_n\}_{n=0}^\infty $ & after (\ref{eqn:def-sp-sol-infty})\\ 
$b_m^n$ & Matrix element for $B_P$: $b_m^n:=\langle B_P e_n,e_m^\Diamond \rangle _{{{\cal H}^\Diamond }}$ & {\bf C2}\\ 
$\ell _0$ & Bandwidth parameter: $b_m^n=0$ for $|m-n|>\ell _0 $ & {\bf C2}\\
$j_0$ & Integer s.t. $b_m^{m+\ell _0 } \ne 0 $
for any integer $m \ge j_0 $ & {\bf C5}\\
$N$ & Dimension of subspace where recursion is executed &{\bf C7} and {\bf Algorithm}\\
$K$ & Dimension of subspace of final approximate solutions & {\bf C7} and {\bf Algorithm}\\
$D$ & Dimension of solution space of $\sum_m b_m^n f_n =0$ & after (\ref{eqn:Haya2})\\
$p_0$ & Integer $p_0:=j_0+\ell _0 -1$& after (\ref{eqn:Haya2})\\
$\|\cdot \|_{\ell ^2,K}$ & `Truncated norm' for number sequences & 
(\ref{lemma:K_norm}) \\ 
\hline
\end{tabular}
\end{center}
\label{tbl:def1}
\end{table}

\subsection{General framework for general linear ordinary differential equation}
\Label{subsec:s21}
In this paper, we consider a 
general linear ordinary differential equation given by the differential operator 
$P(x,{\textstyle \frac{d}{dx}})$ defined on 
the space of $M$-times differentiable functions $C^M(\mathbb{R})$:
\begin{eqnarray}
P(x,{\textstyle \frac{d}{dx}}) f =0, \Label{eqn:C4}
\end{eqnarray}
where $M$ is the order of the differential operator $P(x,{\textstyle \frac{d}{dx}})$.
The main purpose of this paper is to analyze 
the structure of the solution space of (\ref{eqn:C4}) on a given Hilbert space 
${\cal H}$ that densely contains $C^M(\mathbb{R})$.
For this purpose, 
we define the operator $\tilde{A}_P$ as the action of the differential operator $P(x,{\textstyle \frac{d}{dx}})$ with domain
\begin{eqnarray*}
D(\tilde{A}_P):=
\{f\in C^M(\mathbb{R}) \cap {\cal H} | P(x,{\textstyle \frac{d}{dx}}) f \in {\cal H} \}.
\end{eqnarray*}
Then, the linear operator $A_P$ is given as the closed extension of $\tilde{A}_P$ 
with respect to the graph norm ~\cite{ReSi1}.
That is, we treat the structure of the solution space of the differential equation:
\begin{eqnarray*}
A_P f =0.
\end{eqnarray*}
The main goal is to construct an integer-type numerical algorithm for finding non-zero solutions of 
the differential equation given by the differential operator $A_P$
when the original space ${\cal H}$ is contained in a 
larger Hilbert space ${{\cal H}^\Diamond}$ as a set.
Here, the larger Hilbert space ${{\cal H}^\Diamond}$ also densely contains $C^M(\mathbb{R})$.
For this purpose, 
we construct a band-diagonal matrix representation of the differential operator $A_P$
under certain conditions.
In order to obtain a {\it band-diagonal matrix} representation, 
we introduce a linear operator $B_P$ 
from a dense subspace of ${\cal H}$ to ${{\cal H}^\Diamond }$,
which is defined as the closed extension of $\tilde{B}_P$ with respect to the graph norm of the operator $\tilde{B}_P$ defined 
by the action of the differential operator $P(x,{\textstyle \frac{d}{dx}})$ with the following domain:
\begin{eqnarray*}
D(\tilde{B}_P):=
\{f\in C^M (\mathbb{R})\cap {\cal H} | P(x,{\textstyle \frac{d}{dx}}) f \in {{\cal H}^\Diamond } \}.
\end{eqnarray*}

In order to using a band-diagonal structure,
we introduce three conditions
for the quintuplet consisting of
the linear differential operator $P(x,{\textstyle \frac{d}{dx}})$,
the Hilbert spaces ${\cal H}$ and ${{\cal H}^\Diamond }$,
and their CONSs 
$\{e_n \, \}_{n=0}^{\infty}$ and 
$\{e_n^\Diamond \}_{n=0}^{\infty}$,
which is abbreviated to 
$(P(x,{\textstyle \frac{d}{dx}}),{\cal H},
\{e_n \, \}_{n=0}^{\infty},
{{\cal H}^\Diamond },
\{e_n^\Diamond \}_{n=0}^{\infty})$.
These conditions are shown to hold in several examples for $P(x,{\textstyle \frac{d}{dx}})$ later.
In what follows, $\langle \cdot , \, \cdot \rangle _{{{\cal H}^\Diamond }}$
and 
$\langle \cdot ,
 \, \cdot \rangle _{{\cal H}}$ denote the inner products of ${{\cal H}^\Diamond }$ and ${\cal H}$ respectively.
\begin{description}
\item[C1]
For any $n$, $e_n$ belongs to $D(\tilde{B}_P)$.
\item[C2]
There exists an integer $\ell _0 $ 
such that $b_{m}^n:=\langle B_P e_n, e_m^\Diamond \rangle _{{{\cal H}^\Diamond }}=0$ 
when $|n-m|> \ell _0 $.
\item[C3]
There exists a linear operator $C_P$ with domain $D(C_P)$ from a dense subspace of ${{\cal H}^\Diamond }$ to ${\cal H}$
such that $e_m^\Diamond \in D(C_P)$ and 
$\langle B_P f, e_m^\Diamond \rangle _{{{\cal H}^\Diamond }}=\langle f, C_P e_m^\Diamond \rangle _{{\cal H}}$
for $f \in D(\tilde{B}_P)$.
\end{description}

Due to Condition {\bf C3},
the basis $e_m^\Diamond $ belongs to the domain of the adjoint operator $B_P^*$.
Under these conditions, we obtain the following.

\begin{proposition}
\Label{prop:c:matrixeq}
A function $f$ of the kernel of $A_P$ belongs to 
the kernel of $B_P$.
\end{proposition}
This proposition is immediate from the fact that
the domain of $B_P$ includes the domain of $A_P$.

\begin{theorem}
\Label{thm:c:matrixeq}
Assume that the quintuplet $(P(x,{\textstyle \frac{d}{dx}}),{\cal H},\{e_n \, \}_{n=0}^{\infty},
{{\cal H}^\Diamond },\{e_n^\Diamond \}_{n=0}^{\infty})$
satisfies Conditions {\bf C1 - C3}.
For any function $f$ of the kernel of $B_P$
and any $m \in \mathbb{Z}^+ $,
the $\ell ^2$-sequence 
$\{f_n:=\langle f, e_n\rangle_{{\cal H}}\}_{n=0}^\infty $ 
satisfies 
\begin{eqnarray}
\sum_{n=\max(0,\, m-\ell _0)}^{m+\ell _0} b_m^n f_n=0,
\Label{eqn:C1}
\end{eqnarray}
i.e.,
belongs to the linear space
\begin{eqnarray}
V&:=& \bigl \{ \vec{f} :=
\{f_n\}_{n=0}^{\infty}
\,\, \bigl| \, \sum_{n=0}^\infty b_m^n f_n = 0 \,\,\, (m\in\mathbb{Z}^+)\, \bigr\} 
\nonumber \\
&= &\bigl \{ \vec{f} \,\, \bigl| \, \!\!\!\!\! \sum_{n=\max (0, m-\ell _0 )}^{m+\ell _0 }\!\!\!\!\! b_m^n f_n = 0 \,\,\, (m\in\mathbb{Z}^+)\, \bigr\} .
\Label{eqn:def-sp-sol-infty}
\end{eqnarray}
\end{theorem} \par 
\par \noindent The proof of this theorem will be given in Subsection \ref{subsec:s23} of this section.
Due to Condition {\bf C2}, the dimension of $V$ is finite.

However, square-summable number sequences satisfying (\ref{eqn:C1}) do not always correspond to functions in the domain of $A_P$ (hence in the kernel of $A_P$). 
When the linear ordinary differential equation (\ref{eqn:C4}) has 
singular points, i.e., its solution has singular points,
we denote the set of the singular points by $S$.
In this case,
the obtained solutions do not necessarily belong to $C^M(\mathbb{R} )$, 
but they belongs to $C^M(\mathbb{R}\setminus S )$ under some conditions.
So, ${\cal H}$ has to include the space $C^M(\mathbb{R}\setminus S )$.
In order to guarantee that the obtained solutions are true solutions,
we require another condition:
\begin{description}
\item[C4]
For any sequence $\{f_n\}_{n=0}^{\infty} \in 
V \cap \ell ^2(\mathbb{Z}^+) $,
the sum $\displaystyle \sum_{n=0}^N f_n e_n$ converges to a solution $f\in C^M(\mathbb{R}\setminus S )\cap {\cal H}$ of $P(x,\frac{d}{dx})f=0$ 
as $N\to\infty$ for ${\cal H}$-norm. 
\end{description}

Therefore, if the above condition holds,
any a square-summable kernel vector of the band-diagonal matrix $b_m^n$ gives 
a solution of ODE (\ref{eqn:C4}) in 
$C^M(\mathbb{R}\setminus S )\cap {\cal H}$.
In the remainder of this paper, we often use this vector representation instead of a number sequence, for simplicity.


In the non-singular case, the following theorem holds.
\begin{theorem}
\Label{thm:one-to-one}
Assume two assumptions:
(1) the linear ordinary differential equation (\ref{eqn:C4}) has no singular points.
(2) the quintuplet $(P(x,{\textstyle \frac{d}{dx}}),{\cal H},\{e_n \, \}_{n=0}^{\infty},
{{\cal H}^\Diamond },\{e_n^\Diamond \}_{n=0}^{\infty})$
satisfies Conditions {\bf C1}-{\bf C4}.
Then, the map $f\mapsto \{\langle f, e_n\rangle_{{\cal H}}\}_{n=0}^{\infty}$ 
provides a one-to-one correspondence between 
the $\ell ^2$-solutions of (\ref{eqn:C1}) and the solutions in $ C^M(\mathbb{R}) \cap {\cal H} $ of (\ref{eqn:C4}).
\end{theorem} \par 
\par\noindent 
The proof is directly derived from Theorem \ref{thm:c:matrixeq} and Condition {\bf C4} itself as follows.
When a function $f$ in $ C^M(\mathbb{R}) \cap {\cal H} $ satisfies 
the differential equation (\ref{eqn:C4}),
it belongs to $D(\tilde{A}_P) \, (\subset D(A_P) )$ and satisfies $A_P f=0$. 
Then, Theorem \ref{thm:c:matrixeq} and Lemma \ref{prop:c:matrixeq}
guarantee that 
$\{\langle f, e_n\rangle_{{\cal H}}\}_{n=0}^{\infty}$ 
belongs to $V \cap \ell ^2(\mathbb{Z}^+)$.
The reverse argument is immediate from Condition {\bf C4}.

By means of Theorem \ref{thm:one-to-one}, under {\bf C1-C4}, 
the linear differential equation is reduced to 
the simultaneous linear equations (\ref{eqn:C1}) with a `band-diagonal structure' of
bandwidth $2\ell _0+1$. 
That is, under these conditions, 
the problem of finding the solutions in $C^M(\mathbb{R} )\cap {\cal H}$ of the differential equation $P(x,{\textstyle \frac{d}{dx}})f=0$ is equivalent to the problem of finding vectors in the space $V\cap \ell ^2(\mathbb{Z}^+ )$.

Even if the linear ordinary differential equation (\ref{eqn:C4}) has singular points,
we have a modification of Theorem \ref{thm:one-to-one}
if the differential equation (\ref{eqn:C4}) is Fuchsian, whose definition is given as 
follows.

\begin{definition}
An ODE $P(x,\frac{d}{dx})f(x):=
\displaystyle \sum _{m=0}^M p_m(x) (\textstyle \frac{d}{dx})^m f(x)=0$ with analytic functions $p_m(x)$ $(m=0,1,\ldots M)$ is called Fuchsian
if all of its singular points are
regular singular points{\rm ~\cite{CoLe}~\cite{enc}~\cite{Hi}}.
In this case, 
the differential operator 
$P(x,\frac{d}{dx})$ is called Fuchsian.
\end{definition} \par\noindent 
Further, 
in the following, a singular point of ODE $P(x,\frac{d}{dx})f(x)=0$ is
called a singular point of the differential operator $P(x,\frac{d}{dx})$.

As is well known~\cite{enc}~\cite{Hi}, Fuchsian ODEs satisfy the following lemmata:
\begin{lemma}\Label{lem:1}
Assume that a Fuchsian ODE $\displaystyle \sum _{m=0}^M p_m(x) (\textstyle \frac{d}{dx})^m f(x)=0$ with analytic functions $p_m(x)$ $(m=0,1,\ldots M)$
has singular points $z_1, z_2, \ldots z_{N}$.
Then, 
the functions $\displaystyle \prod _{n=1}^{N} (x-z_n)^{M-m}\, \frac{p_m(x)}{p_M(x)}$ $(m=0,1,\ldots M-1)$ 
are holomorphic at $z_1, z_2, \ldots z_{N_1}$.
\end{lemma}

\begin{lemma}\Label{lem:2}
Assume that a Fuchsian ODE $\displaystyle \sum _{m=0}^M p_m(x) (\textstyle \frac{d}{dx})^m f(x)=0$ has 
holomorphic coefficient functions $p_m(x)$ $(m=0,1,\ldots M)$ on $\mathbb{R}$.
Then, 
the set $S$ of its singular points
is given by $p_M^{-1}(0)$.
\end{lemma}

As a corollary, we obtain the following:

\begin{corollary}\Label{cor:1}
Assume that a differential operator  $
Q(x,\frac{d}{dx}):=\displaystyle \sum _{m=0}^M q_m(x) (\textstyle \frac{d}{dx})^m $ with analytic functions $q_m(x)$ $(m=0,1,\ldots M)$
has singular points $z_1, z_2, \ldots z_{N}$
and  
all of zero points of the coefficient function $q_M(x)$ are
$z_1, z_2, \ldots z_{N}$ whose multiplicity are 
not smaller than $M$. Moreover, assume that the ODE $Q(x,\frac{d}{dx})f=0$ is Fuchsian. 
Then,
there exist holomorphic functions $\tilde{q}_m(x)$ $(m=0,1,\ldots M-1)$ on $\mathbb{R}$
such that $\displaystyle q_m(x)=\tilde{q}_m(x)\prod _{n=0}^{N_1} (x-z_n)^{m}$ 
and the set $S$ of its singular points
is given by $q_M^{-1}(0)$.
\end{corollary}

We additionally assume the conditions:
\begin{description}
\item[C1$^+$]
There exists a positive function $\upsilon$ in $C^M(\mathbb{R} )$ s.t. $\displaystyle \langle f,g \rangle _{\cal H}=\int_{-\infty }^\infty \!\!\! f(x)\overline{g(x)}\upsilon(x)dx$.
\item[C2$^+$]
There exists a positive function ${\upsilon^\Diamond}$ in $C^M(\mathbb{R} )$ s.t. $\displaystyle \langle f,g \rangle _{{\cal H}^\Diamond}
=\int_{-\infty }^\infty \!\!\! f(x)\overline{g(x)}{\upsilon^\Diamond}(x)dx$.
\end{description}
When Condition {\bf C1$^+$} holds, 
${\cal H}$ always includes the space $C^M(\mathbb{R}\setminus S )$
because the set $S$ has zero measure.
Now, Theorem \ref{thm:one-to-one} can be replaced by the following theorem:

\begin{theorem}
\Label{thm:one-to-one2}
Let $Q(x,{\textstyle \frac{d}{dx}})$ be the Fuchsian differential operator 
satisfying the conditions of Corollary \ref{cor:1}.
When the quintuplet $(Q(x,{\textstyle \frac{d}{dx}}),
{\cal H},\{e_n \, \}_{n=0}^{\infty},
{{\cal H}^\Diamond },\{e_n^\Diamond \}_{n=0}^{\infty})$
satisfies Conditions {\bf C1}-{\bf C4}, {\bf C1$^+$}, and {\bf C2$^+$}, 
the map $f\mapsto \{\langle f, e_n\rangle_{{\cal H}}\}_{n=0}^{\infty}$ 
provides a one-to-one correspondence between 
the $\ell ^2$-solutions of (\ref{eqn:C1}) 
and the solutions in 
$C^M(\mathbb{R}\setminus S) \cap {\cal H} $ of (\ref{eqn:C4})
concerning $B_Q$.
\end{theorem} \par 

Hence, 
when the given conditions hold,
in order to solve the Fuchsian linear ODE,
it is sufficient to extract the subspace $V \cap \ell ^2(\mathbb{Z}^+)$
as well as in the non-singular case.

However, a general Fuchsian differential operator 
$P(x,{\textstyle \frac{d}{dx}}):=\displaystyle \sum _{m=0}^M p_m(x) (\textstyle \frac{d}{dx})^m$
does not necessarily satisfy the above condition.
In this case, Theorem \ref{thm:one-to-one2} can be applied in the following way.
Assume that 
all coefficient functions $p_m(x)$ $(m=0,1,\ldots M)$ are holomorphic
and 
$p_M(x)$ has zero points $z_1, \ldots, z_N \in S$ with 
the multiplicity $\mu_1, \ldots, \mu_N$, respectively.
Then, Lemma \ref{lem:1} guarantees the inequality $\mu_n \le M$.
Then, the 
differential operator
\begin{eqnarray}
Q(x,\frac{d}{dx}):=\displaystyle \prod_n (x-z_n)^{M-\mu_n}P(x,\frac{d}{dx})
\Label{eq:4-5-1}
\end{eqnarray}
is Fuchsian and satifies the conditions of Corollary \ref{cor:1}.
So, we can apply Theorem \ref{thm:one-to-one2} to the differential operator $Q(x,\frac{d}{dx})$ 
in stead of $P(x,\frac{d}{dx})$.

Since Condition {\bf C4} is assumed in Theorem \ref{thm:one-to-one2},
it is sufficient to show the following theorem, which will be shown in Subsection \ref{subsec:s24}.
That is, the combination of {\bf C4} and Theorems \ref{thm:c:matrixeq} and \ref{thm:s35}
yields Theorem \ref{thm:one-to-one2}.

\begin{theorem}\Label{thm:s35}
Let $Q(x,{\textstyle \frac{d}{dx}})$ be the Fuchsian differential operator 
satisfying the conditions of Corollary \ref{cor:1}.
Assume that the quintuplet $(Q(x,{\textstyle \frac{d}{dx}}),
{\cal H},\{e_n \, \}_{n=0}^{\infty},
{{\cal H}^\Diamond },\{e_n^\Diamond \}_{n=0}^{\infty})$
satisfies Conditions {\bf C1}-{\bf C3}, {\bf C1$^+$}, and {\bf C2$^+$}.
Then, the solutions in 
$C^M(\mathbb{R}\setminus S) \cap {\cal H} $ of (\ref{eqn:C4}) always
give functions of the kernel of $B_Q$.
\end{theorem}

In this theorem, 
the condition for the multiplicity of zero points is crucial
because it is needed to expand the domain $D(B_Q)$ sufficiently.

\begin{remark}
\label{rem:nonsym}
In general, the inner product $\langle \cdot , \, \cdot \rangle _{{\cal H}}$ of ${\cal H}$ does not coincide with the restriction 
on ${\cal H}$ of the inner product $\langle \cdot , \, \cdot \rangle _{{{\cal H}^\Diamond }}$ of ${{\cal H}^\Diamond }$.
Further, 
$e_n^\Diamond $ does not necessarily belong to the domain $D(A_P^*)$.
In order to characterize Condition {\bf C3},
we consider the special case when 
${{\cal H}^\Diamond }={\cal H}$, $e_n=e_n^\Diamond $, and
$\langle A_P e_n , e_m\rangle _{{{\cal H}^\Diamond }} =\overline{\langle A_P e_m , e_n\rangle _{{{\cal H}^\Diamond }}}$.
Note that the condition ${{\cal H}^\Diamond }={\cal H}$ implies that $A_P=B_P$.
In this case, 
if the operator $A_P$ is symmetric,
Condition {\bf C3} holds.
In other words, if the operator $A_P$ is not symmetric,
Condition {\bf C3} does not necessarily hold.
In such a case, 
if we define another linear operator $\tilde{A}_P'$ whose domain is the linear expansion of 
$\{e_n \}$, its closed extension $A_P'$ is symmetric, 
and the solution function of $A_P'f=0$ satisfies the simultaneous linear equations (\ref{eqn:C1}).
That is,
if a general non-zero solution function of $A_Pf=0$ does not belong to the domain of $A_P'$, 
this solution does not necessarily satisfy (\ref{eqn:C1}). Later, in Remark \ref{rem:nonsym_ex}, we give an example of the latter case.
\end{remark}

The remaining tasks are divided into two parts:
The first part concerns the general theory for our algorithm for solving a linear ordinary differential equation based on several conditions.
This part is called general theory part.
The second part concerns how to apply the above general theory for 
several wide classes of linear ordinary differential equations.
This part is called application part.

{\bf General part}
\begin{description}
\item[Task 1]
(Subsection \ref{subsec:s22})
Giving a recursive algorithm for band-diagonal-type simultaneous linear equations
under Conditions {\bf C1}-{\bf C4}.
This algorithm requires additional condition {\bf C5},
which will be given in Subsection \ref{subsec:s22}.
The additional explanation for this algorithm is given in Subsection \ref{ss23}.

\item[Task 2]
(Subsection \ref{ss23})
Giving an integer-type algorithm realizing the above algorithm
by adding Condition {\bf C6}.

\item[Task 3]
(Subsection \ref{subsec:s23})
Proof of Theorem \ref{thm:c:matrixeq}:
Showing that any $C^M(\mathbb{R})$ solution corresponds to an element of $V$.

\item[Task 4]
(Subsection \ref{subsec:s24})
Proof of Theorem \ref{thm:s35}:
Showing the one-to-one correspondence 
between 
$V \cap \ell ^2(\mathbb{Z}^+)$ and $C^M(\mathbb{R}\setminus S)\cap {\cal H}$
with the Fuchsian differential operator $Q$ given in (\ref{eq:4-5-1}).

\item[Task 5]
(Subsections \ref{ss43})
Showing the convergence of the algorithm given in Task 1.

\item[Task 6]
(Subsections \ref{ss44})
Showing that all of $\ell ^2$ components in $V$ can be extracted by 
the algorithm given in Task 1 with additional condition.

\end{description}

{\bf Application part}
\begin{description}
\item[Task 7]
(Subsections \ref{subsec:s41}-\ref{subsec:s44})
Constructing
${\cal H}$, ${{\cal H}^\Diamond }$, and their CONSs satisfying Conditions  
{\bf C1}-{\bf C6}, {\bf C1$^+$}, and {\bf C2$^+$}
for a differential operator $P(x,\frac{d}{dx})$ with polynomial coefficient functions $p_m(x)$.

\item[Task 8]
(Subsection \ref{subsec:s41-1})
Explaining how to apply the above method to a differential operator $R(x,\frac{d}{dx})$
with rational coefficient functions $r_m(x)$.


\end{description}

In order to apply our method to a Fuchsian linear ODE $P(x,\frac{d}{dx})f=0$ 
with polynomial coefficient functions,
it is enough to construct
${\cal H}$, ${{\cal H}^\Diamond }$, and their CONSs satisfying Conditions  
{\bf C1}-{\bf C6}, {\bf C1$^+$}, and {\bf C2$^+$}
with the Fuchsian linear differential operator $Q(x,\frac{d}{dx})$ given in (\ref{eq:4-5-1}).
Since the Fuchsian linear differential operator $Q(x,\frac{d}{dx})$
has polynomial coefficient functions,
we can apply 
{\bf Task 7} with replacing $P(x,\frac{d}{dx})$ by $Q(x,\frac{d}{dx})$.
Due to Theorem \ref{thm:one-to-one2},
any solution of $P(x,\frac{d}{dx})f(x)=0$ in $C^M(\mathbb{R}\setminus p_M^{-1}(0))
\cap {\cal H}$ can be obtained by this method.

\subsection{Recursive algorithm for band-diagonal-type simultaneous linear equations}
\Label{subsec:s22}

In the next step,
we consider the algorithm for $\ell ^2$-solution of the band-diagonal simultaneous linear equations (\ref{eqn:C1}).
In this subsection, we briefly describe 
the structure of our algorithm for this problem
and 
explain how to avoid the usual pitfalls of this method. 

From {\bf C3}, the simultaneous linear equations (\ref{eqn:C1}) have a `band-diagonal structure' with bandwidth $2\ell _0+1$. This type of system of simultaneous linear equations can be solved easily.
The simultaneous linear equations (\ref{eqn:C1}) with {\bf C2} have at least $\ell _0 $ linearly independent algebraic solutions.
The linearly independent solutions of the solution space $V$ defined in (\ref{eqn:def-sp-sol-infty}) 
can be solved recursively when the following condition holds
for the quintuplet $(P(x,{\textstyle \frac{d}{dx}}),{\cal H},\{e_n \, \}_{n=0}^{\infty},
{{\cal H}^\Diamond },\{e_n^\Diamond \}_{n=0}^{\infty})$.
\begin{description}
\item[C5]
There exists an integer $j_0\in \mathbb{Z}^+$ such that $b_m^{m+\ell _0 } \ne 0 $
for any integer $m \ge j_0 $ $ (m\in \mathbb{Z} ^+)$.
\end{description}
In this case, the dimension $D$ of $V$ is equal to that of 
\begin{eqnarray}
\Pi_{p_0} V=\bigl \{ \{f_n\}_{n=0}^{p_0} \, \bigl| \, 
\sum_{n=0}^{p_0} b_m^n f_n = 0 \,\, (m=0,1,...,j_0-1)\, \bigr\},\Label{eqn:Haya2}
\end{eqnarray}
where $p_0:=j_0+\ell _0 -1$ and 
the truncation operator $\Pi_m$ is defined by 
\begin{eqnarray}
(\Pi_m \vec{f})_n = \left\{\begin{array}{@{\,}ll} f_n & \,\,\,\, 
(n\le m) \\ 0 & \,\,\,\, (n>m)\,\, . \end{array}\right.
\Label{eqn:1}
\end{eqnarray}
In the following, for simplicity, we sometimes identify $\Pi_m\vec{f}$ with the corresponding  $(m+1)$-dimensional vector.  
This band-diagonal matrix $b_{m}^n$ is illustrated by Figure \ref{f1}.

\begin{figure}[htbp]
\begin{center}
\scalebox{1.0}{\includegraphics[scale=0.4]{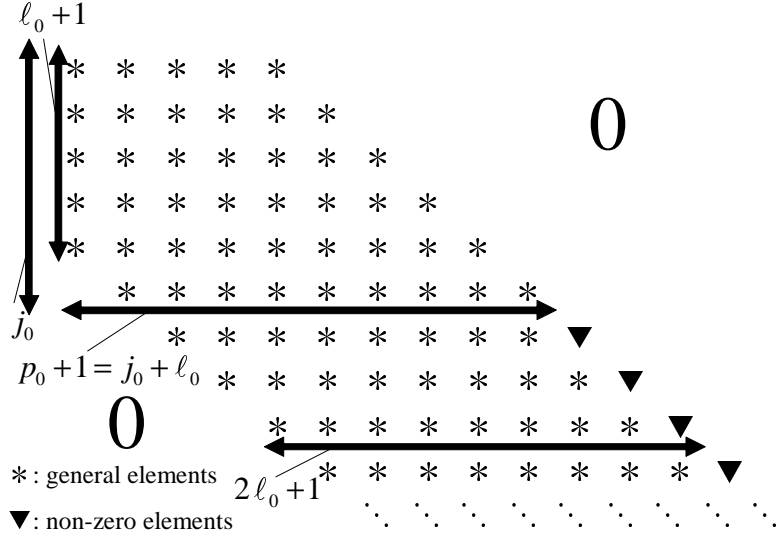}}
\end{center}
\caption{Figure of band-diagonal matrix $b_{m}^n$}
\Label{f1}
\end{figure}%

So, 
we define $D$ linearly independent sequences $\vec{F}^{(1)}=\{F^{(1)}_n\}_{n=0}^{\infty}, 
\ldots, \vec{F}^{(D)}=\{F^{(D)}_n\}_{n=0}^{\infty}$ by the following procedure:
the first ${p_0}$ elements of all sequences $\vec{F}^{(1)}, \ldots, \vec{F}^{(D)}$ 
by $D$ linearly independent vectors of $V$.
The remaining elements $F_n^{(d)}$ with $n\ge {p_0}+1$ are calculated by the recursion 
\begin{eqnarray}
F_n^{(d)} = -\, \frac{1}{b_{n-\ell _0 }^n} \sum_{m=n-2\ell _0 }^{n-1} b_{n-\ell _0}^m F_m^{(d)} \,\, , 
\Label{eqn:recursive-alg} \end{eqnarray}
because $b_{n-\ell _0 }^n\ne 0$ there. 
(The first procedure to find a basis system is easy; it is to solve a system of finite-dimensional simultaneous linear equations.)
The following theorem follows directly from the construction of 
$\vec{F}^{(1)}, \ldots, \vec{F}^{(D)}$. 
Therefore, we obtain the following proposition:
\begin{proposition}
\Label{thm:recursive-alg}
Under {\bf C5}, 
the algebraic solution space $U $ can be spanned by $\vec{F}^{(1)}, \ldots, \vec{F}^{(D)}$.
That is, any algebraic solution $\vec{f}$ of (\ref{eqn:C1})
can be obtained by a linear combination of the basis sequences 
$\vec{F}^{(1)}, \ldots, \vec{F}^{(D)}$.
\end{proposition} \par 

However, here is an important pitfall. From the existence and uniqueness theorems, 
there are $M$ linearly independent true solutions in $C^M(\mathbb{R})$
when there is no singular point, i.e., $S=\emptyset$. 
In this case, therefore, in $C^M(\mathbb{R})\cap {\cal H}$, the number of linearly independent solutions is not greater than $M$, 
which is smaller than $D$
in a typical example given in Section \ref{sec:CSI}. 
Even though there exist singular points, 
due to Theorem \ref{thm:c:matrixeq},
the solutions in $\ell ^2(\mathbb{Z}^+)$ of the simultaneous linear equations correspond to the true solutions in $C^M(\mathbb{R}\setminus S)\cap {\cal H}$ of the ODE. 
However, when a solution of the simultaneous linear equations does not belong to $\ell ^2(\mathbb{Z}^+)$, 
it usually does not correspond to a true solution in $C^M(\mathbb{R}\setminus S)\cap {\cal H}$ of the ODE. 
Therefore, we have to be careful to this differentiate.

In general, the solution $\vec{f}$ obtained by the above recursion 
is a linear combination of these three kinds of components, and 
it is not so easy to extract 
the component corresponding to the true solution in 
$C^M(\mathbb{R}\setminus S)\cap {\cal H}$.
In the following, we propose a method to extract the $\ell ^2$-component, i.e.,
one element of the subset $V \cap \ell ^2(\mathbb{Z}^+) $. 
For these purposes, 
we choose a bounded bilinear form $\Omega(\vec{f}, \, \vec{g})$ on 
$\ell ^2(\mathbb{Z}^+)\times \ell ^2(\mathbb{Z}^+)$ 
(and the corresponding quadratic form $\Omega(\vec{f}):=\Omega(\vec{f}, \, \vec{f})$ on 
$\ell ^2(\mathbb{Z}^+)$ ) and the integers $K$ and $N$ 
satisfying
\begin{eqnarray}
& ^\forall \vec{f} \in \ell ^2(\mathbb{Z}^+), \quad 
\Omega(\vec{f})\ge \|\vec{f}\|_{\ell ^2}^2 :=
\sum_{n=0}^{\infty} |f_n|^2,\quad
\Label{OC6} \\
& N \ge K \ge j_0+\ell _0 -1,
\Label{OC7}
\end{eqnarray} 
and define the ratio and its minimum:
\begin{eqnarray}
\sigma_{K,N}^{(\Omega )}(\vec{f}) 
&:=&
\frac{\Omega (\Pi_N \vec{f})}
{\|\vec{f}\|_{\ell ^2,K}^2} \quad \hbox{ for } \vec{f}  \in V\setminus\{0\}\\
\underline{\sigma_{K,N}^{(\Omega )}}
&:=&
\min_{\vec{f} \in V\setminus\{0\}} \sigma_{K,N}^{(\Omega )}(\vec{f}) .
\Label{eqn:def_ratio}
\end{eqnarray}
This definition is well defined because 
Conditions {\bf C2}, {\bf C5} and (\ref{OC7}) guarantee the relation
\begin{eqnarray}
\| \vec{f}\|_{\ell ^2,K}:=\|\Pi_K\vec{f}\|_{\ell ^2}
>0
\quad \hbox{ for } \vec{f}  \in V\setminus\{0\}.
\Label{lemma:K_norm}
\end{eqnarray}
Similarly, 
we define 
\begin{eqnarray}
\sigma_{K,\infty}^{(\Omega )}(\vec{f}) 
&:=& 
\frac{\Omega (\vec{f})}
{\|\vec{f}\|_{\ell ^2,K}^2} \quad \hbox{ for } \vec{f}  \in V\setminus\{0\}\\
\underline{\sigma_{K,\infty}^{(\Omega )}}
&:=&
\min_{\vec{f} \in V\setminus\{0\}} \sigma_{K,\infty}^{(\Omega )}(\vec{f}) .
\end{eqnarray}
Hence, 
our solution space 
$V_{K}$ is given with the following condition
\begin{eqnarray}
V_{K}\subset 
\Pi_{K}
\left((\sigma_{K,N}^{(\Omega )})^{-1}[0,c \underline{\sigma_{K,N}^{(\Omega )}} ]\right) 
\cup \{0\}
\Label{4-3-1}
\end{eqnarray}

Now, we introduce our algorithm to approximately obtain
the $\ell ^2(\mathbb{Z}^+) $ components of $V$:
\par\noindent 
\vspace{2mm}
\par\noindent 
{\bf Algorithm}
\vspace{3mm}
$\,\,\,\,\,\,\,\, $ 
\begin{description}
\item[Step 1]
Calculation of basis vectors of $\Pi_{p_0} V$:
\item[] 
\hspace{5mm} Find a basis system $\{F_n^{(1)}\}_{n=0}^{p_0} ,\ldots, \{F_n^{(D)}\}_{n=0}^{p_0}$
 for $\Pi_{p_0} V$ in (\ref{eqn:Haya2}) 
by Gaussian elimination, where $D$ is determined by its result. 
This is easy because $p_0$ is small.
\item[Step 2] Recursive calculation of basis vectors of $\Pi_{n} V$ $(p_0+1\le n \le N)$: 
\item[]
\hspace{5mm} Iterate the recursion (\ref{eqn:recursive-alg}) for $n=p_0+1,p_0+2,\ldots , N$, in order to obtain a basis system $\{F_n^{(1)}\}_{n=0}^{N} ,\ldots, \{F_n^{(D)}\}_{n=0}^{N}$ for $\Pi_N V$.
\item[Step 3] Removal of components from $\Pi_K V$ 
corresponding to non-$\ell ^2$-components in $V$: 
\item[]
\hspace{5mm} Find a linear subspace 
$V_{K}$
of 
$\Pi_{K}\left((\sigma_{K,N}^{(\Omega )})^{-1}[0,c \underline{\sigma_{K,N}^{(\Omega )}} ]\right) \cup \{0\}$.
This process can be done as follows.
\begin{description}
\item[Step 3.1]
Find a non-zero vector $\{G_n^{(1)}\}_{n=0}^{K}$ in 
$\Pi_{K}\left((\sigma_{K,N}^{(\Omega )})^{-1}[0,c \underline{\sigma_{K,N}^{(\Omega )}} ]\right) $.
\item[Step 3.2]
Find a non-zero vector $\{G_n^{(2)}\}_{n=0}^{K}$ in 
$\Pi_{K}\left((\sigma_{K,N}^{(\Omega )})^{-1}[0,c \underline{\sigma_{K,N}^{(\Omega )}} ]\right) 
\cap <\{G_n^{(1)}\}_{n=0}^{K}>^{\perp,K}$.
\item[]$\vdots$
\item[Step 3.$n$]
Find a non-zero vector $\{G_n^{(i)}\}_{n=0}^{K}$ in 
$\Pi_{K}\left((\sigma_{K,N}^{(\Omega )})^{-1}[0,c \underline{\sigma_{K,N}^{(\Omega )}} ]\right) 
\cap <\{G_n^{(1)}\}_{n=0}^{K},\ldots, \{G_n^{(i-1)}\}_{n=0}^{K}>^{\perp,K}$.
\end{description}
When 
$\Pi_{K}\left((\sigma_{K,N}^{(\Omega )})^{-1}[0,c \underline{\sigma_{K,N}^{(\Omega )}} ]\right) 
\cap <\{G_n^{(1)}\}_{n=0}^{K},\ldots, \{G_n^{(i-1)}\}_{n=0}^{K}>^{\perp,K}$
is an empty set, we stop this process.
Here,
$~^{\perp,K}$ denotes the orthogonal space concerning the inner product of
$\ell^2(\{0,\ldots, m\})$.
\end{description}
\vspace{3mm}

When the purpose is to
 calculate the truncated elements $\ell ^2$ solution $\Pi_K(V \cap \ell ^2(\mathbb{Z}^+))$,
the error is evaluated by the norm concerning inner product 
$\langle \vec{x},\vec{y}\rangle_{\ell ^2,K}:=
\langle \Pi_K\vec{x},\Pi_K\vec{y}\rangle_{\ell ^2}$.
Denoting the 
the projection to 
$W$ concerning this inner product by $P_{W,K}$,
we can evaluate the accuracy of our result $V_{K}$ by
\begin{eqnarray}
\sup_{\vec{x} \in V_{K}\setminus\{0\}}
\frac{ \|P_{V\cap \ell ^2(\mathbb{Z}^+),K} \vec{x} - \vec{x}\|_{\ell ^2,K} }{\|\vec{x}\|_{\ell ^2,K} }.
\end{eqnarray}
However, 
the subspace $V_{K}$ is not uniquely defined, and is chosen with the condition
(\ref{4-3-1}).
The accuracy of our algorithm should be evaluated with the worst case as follows:
\begin{eqnarray}
\nonumber 
&& \sup_{V_{K}\subset 
\Pi_{K}
((\sigma_{K,N}^{(\Omega )})^{-1}[0,c \underline{\sigma_{K,N}^{(\Omega )}} ]) \cup \{0\}} \,\, 
\sup_{\vec{x} \in V_{K}\setminus\{0\}}
\frac{ \|P_{V,K} \vec{x} -\vec{x}\|_{\ell ^2,K} }{\|\vec{x}\|_{\ell ^2,K} }\\
\nonumber 
&=& 
\sup_{
\vec{x} \in 
\Pi_{K}
((\sigma_{K,N}^{(\Omega )})^{-1} [0,c \underline{\sigma_{K,N}^{(\Omega )}} ]) }
\frac{ \|P_{V,K} \vec{x} -\vec{x}\|_{\ell ^2,K} }{\|\vec{x}\|_{\ell ^2,K} }\\
&=& 
\sup_{
\vec{x} \in 
(\sigma_{K,N}^{(\Omega )})^{-1} [0,c \underline{\sigma_{K,N}^{(\Omega )}} ] }
\frac{ \|P_{V,K} \vec{x} -\vec{x}\|_{\ell ^2,K} }{\|\vec{x}\|_{\ell ^2,K} }.
\end{eqnarray}
It can be shown that this value goes to $0$ as follows.

\begin{theorem}\Label{thm:conv1}
For fixed $K$, when $N$ goes to infinity,
the convergence
\begin{eqnarray}
\sup_{
\vec{x} \in 
(\sigma_{K,N}^{(\Omega )})^{-1} [0,c \underline{\sigma_{K,N}^{(\Omega )}} ] }
\frac{ \|P_{V,K} \vec{x} -\vec{x}\|_{\ell ^2,K} }{\|\vec{x}\|_{\ell ^2,K} }
\to 0
\end{eqnarray}
holds.
\end{theorem}

Hence, 
since $\Pi_K \vec{x} $
is close to $\vec{x}$ for $\vec{x} \in V \cap \ell ^2(\mathbb{Z}^+)$,
the above theorem guarantees that
our algorithm gives a subspace of $V$ whose elements are close to elements of 
$V \cap \ell ^2(\mathbb{Z}^+)$.

Conversely, 
the above theorem cannot guarantee that
all of elements is approximated by elements of $V_K$.
For this purpose, we have to guarantee that 
$P_{V,K}V_K $ is equal to $\Pi_K( V \cap \ell ^2(\mathbb{Z}^+))$,
which is shown as follows.

\begin{theorem}\Label{thm:conv2}
Assume that $\Omega(\vec{x})=\|\vec{x}\|_{\ell ^2}$.
When we choose sufficiently large numbers $N_0$ and $c_0$,
then for any $N \ge N_0$ and $c\ge c_0$, we have 
\begin{eqnarray}
P_{V,K}V_K = \Pi_K( V \cap \ell ^2(\mathbb{Z}^+))\Label{eq:2}
\end{eqnarray}
for any choice of $V_K$.
\end{theorem}

Therefore, when we choose $\Omega(\vec{x})=\|\vec{x}\|_{\ell ^2}$
and sufficiently large numbers $N$ and $c$,
we obtain a subspace $V_K \subset V$ close to $\Pi_K (V \cap \ell ^2(\mathbb{Z}^+))$.

>From the above reason, the above choice of $\Omega$ suits our purpose.
However, 
when we choose the quadratic form $\Omega$ different from $\|\vec{x}\|_{\ell ^2}$,
the convergence is improved in several specific example.
That is,
for sake of practical accuracy of the solutions obtained, 
a condition concerning a kind of sensitivity of 
the quadratic form $\Omega$ will be required later.
The details will be given in our paper ~\cite{paper3}, 
and we will use such a sensitive quadratic form in the numerical examples in Section \ref{sec:nm}.

\subsection{Realization by an integer-type algorithm}
\Label{ss23}

When the quintuplet $(P(x,{\textstyle \frac{d}{dx}}),{\cal H},\{e_n \, \}_{n=0}^{\infty},
{{\cal H}^\Diamond },\{e_n^\Diamond \}_{n=0}^{\infty})$.
satisfies the following condition,
the above algorithm can be realized 
by the following integer-type algorithm
with a small modification.
\begin{description}
\item[C6]
There exists a complex number $\gamma \in \mathbb{C}$ such that
$\gamma \, b_{m}^n\in \mathbb{Q} + \mathbb{Q} i, \,\,\, (m,n\in \mathbb{Z} )$.
\end{description}

The most crucial part is {\bf Step 3}.
To execute {\bf Step 3}, a simple method is to calculate a vector $\vec{f}$ which minimizes the ratio $\sigma_{K}^{(\Omega)}(\vec{f})$. Usually, with the matrices $A$ and $B$ defined by 
$(A)_{i\, j}:=(\Pi_K\vec{F}^{(i)},\, \Pi_K\vec{F}^{(j)})_{\ell ^2}$ and $(B)_{i\, j}:=\Omega(\Pi_N\vec{F}^{(i)}, \Pi_N\vec{F}^{(j)})$, this minimization can be performed exactly by calculating an eigenvector of the matrix $A^{-\frac{1}{2}}BA^{-\frac{1}{2}}$ associated with the minimum eigenvalue. However, this normal method is quite difficult to apply, 
because the matrices $A$ and $B$ are usually very close to a singular matrix with rank $1$ due to the most diverging components in $V$, 
and hence this usual method is particularly subject to 
the `canceling' due to round-off errors.  
In order to avoid this, many calculations are required, 
if we try to find the optimal vector with high 
accuracy by the usual methods.

We can avoid so many calculations and 
the `canceling' due to round-off errors
in this minimization by carrying out the following quasi-minimization.
As is proved by a geometrical discussion of the convex set in a more general framework in~\cite{paper3}, by means of the Schwarz inequality,
an arbitrary orthogonal basis system for $\Pi_N V$ 
with respect to the inner product $\langle \cdot, \, \cdot \rangle _\Omega$ 
contains at least one vector $\vec{f}$ such that 
$\vec{f} \in 
\Pi_{K}\left((\sigma_{K,N}^{(\Omega )})^{-1}[0,c \underline{\sigma_{K,N}^{(\Omega )}} ]\right) $.
Hence, we have only to take the basis vector with the minimum ratio 
$\sigma_{K,N}^{(\Omega)}(\vec{f})$
among this basis system. 

However, 
in order to take an orthogonal basis system, 
we need exact orthogonalization, which requires a large amount of calculations. 
To avoid this problem, we propose an alternative method based on a kind of integer-type quasi-orthogonalization of a $D$-dimensional `lattice', where the angles between the final basis vectors are not distant by more than $\zeta $ $(<<1)$ from being exactly orthogonal, by which 
$\vec{f} \in 
\Pi_{K}\left((\sigma_{K,N}^{(\Omega )})^{-1}[0,c \underline{\sigma_{K,N}^{(\Omega )}} ]\right) $
is guaranteed for the basis vector with minimum ratio $\sigma_{K}^{(\Omega)}(\vec{f})$. 
This method is somewhat similar to the `lattice reduction problem'~\cite{lat1}~\cite{lat2}, which is well known as an NP-hard problem if we require exact minimization of the lattice. 
However, our alternative method aims at closeness to orthogonality rather than exact minimization of the lattice, only with 
a small amount of  
calculations, by means of a quasi-orthogonalization algorithm which does not increase the integers used for the numerators and the common denominator of complex rational numbers except for special cases with bad final orthogonality~\cite{paper3}. 
So, we can perform {\bf Step 3.1} with few calculations.

When the dimension $D_{\ell ^2}$ of $V \cap \ell ^2(\mathbb{Z}^+)$ is strictly greater than $1$,
we need to perform {\bf Step 3.2}, $\ldots$, {\bf Step 3.$D_{\ell ^2}$}.
In the first process in these steps, 
we need 
to calculate $\Pi_{K}\left((\sigma_{K,N}^{(\Omega )})^{-1}[0,c \underline{\sigma_{K,N}^{(\Omega )}} ]\right) 
\cap <\{G_n^{(1)}\}_{n=0}^{K},\ldots, \{G_n^{(i-1)}\}_{n=0}^{K}>^{\perp,K}$,
which requires 
orthogonalization concerning the inner product $\langle ~,~\rangle_{\ell ^2, K}$.
We apply the above quasi-orthogonalization algorithm to this orthogonalization.
If the second vector $\{G_n^{(2)}\}_{n=0}^{K}$ 
is not orthogonal to the first vector $\{G_n^{(1)}\}_{n=0}^{K}$,
linear combinations of 
$\{G_n^{(1)}\}_{n=0}^{K}$ and $\{G_n^{(2)}\}_{n=0}^{K}$
may not belong to 
$\Pi_{K}\left((\sigma_{K,N}^{(\Omega )})^{-1}[0,c \underline{\sigma_{K,N}^{(\Omega )}} ]\right) $
However, as is shown Section 6 of \cite{paper3},
if $\{G_n^{(2)}\}_{n=0}^{K}$ 
is sufficiently close to a vector orthogonal to the first vector $\{G_n^{(1)}\}_{n=0}^{K}$,
any linear combination of both vector belongs to
$\Pi_{K}\left((\sigma_{K,N}^{(\Omega )})^{-1}[0,c \underline{\sigma_{K,N}^{(\Omega )}} ]\right) $.
Repeating this procedure up to
the $D_{\ell ^2}$ times,
we can find a linear subspace $V_K$ as a subset of 
$\Pi_{K}\left((\sigma_{K,N}^{(\Omega )})^{-1}[0,c \underline{\sigma_{K,N}^{(\Omega )}} ]\right) $.

Therefore, we can realize the above algorithm by an integer-type algorithm with few calculations.

\subsection{Possibility of the estimation of accuracy}
In numerical methods, it is important whether or not we can determine the accuracy of numerical results. For our method, we will give an upper bound of the norm of total errors in~\cite{paper3}. This error bound is a function only of the norm of the truncation error due to the components outside the subspace ${\cal H}^{(K)}={\rm Span}(e_0,e_1,\ldots ,e_K)$, and all the other parameters for the bound than this truncation error can be calculated using only the numerical results without requiring any knowledge of the true solutions.

\section{Function spaces and basis systems used in our method}
\Label{sec:CSI}

In this section, 
firstly, we explain what spaces are used for ${\cal H}$ and ${{\cal H}^\Diamond }$ as well as what basis function systems are used in our algorithm 
in Subsection \ref{subsec:s41}
when 
$P(x,{\textstyle \frac{d}{dx}})=
\sum_{m=0}^{M}p_m(x)(\frac{d}{dx})^m$ 
can be written as a polynomial in $x$ and ${\textstyle \frac{d}{dx}}$.
So, in the singular case,
the set $S$ of singular points is given by $p_M^{-1}(0)$, that is,
it is equal to the set of zero points of $p_M$.
Next, 
we explain that the presented examples satisfy 
Conditions {\bf C1}-{\bf C5}, {\bf C1$^+$}, and {\bf C2$^+$}
in Subsections \ref{subsec:s41} (for {\bf C1}, {\bf C1$^+$}, and {\bf C2$^+$}),
\ref{subsec:s42} (for {\bf C2}, {\bf C5}, and {\bf C6}), 
\ref{subsec:s43} (for {\bf C3}), 
and \ref{subsec:s44} (for {\bf C4}). 
In Subsection \ref{subsec:s41-1}, we explain how to apply our method 
to the non-polynomial case.

\begin{table}[bht]
\caption{Notation in Section 4}
\begin{center}
\begin{tabular}{|l@{}|l@{}|@{}c@{}|}
\hline 
$L_{(k_0)}^2 (\mathbb{R})$ & Input space (concrete choice) & (\ref{eqn:sp}) and (\ref{eqn:spaces_CSI1})\\
$L_{({k_0^\Diamond })}^2 (\mathbb{R})$ & Output space (concrete choice)& (\ref{eqn:sp}) and (\ref{eqn:spaces_CSI2}) \\
$k_0$ & Integer parameter for input space & (\ref{eqn:sp}) and (\ref{eqn:spaces_CSI1})\\
${k_0^\Diamond }$ & Integer parameter for output space & (\ref{eqn:sp}) and (\ref{eqn:spaces_CSI2})\\
$s_0$& Obligatory minimum of difference $k_0-k_0^\Diamond $ & after (\ref{eqn:spaces_CSI2}) \\
$\psi_{k,\, {\ddot{n}}}(x)$ & Wavepacket function used for bases& (\ref{eqn:def-psi-CSI})\\
$\ddot{n}_{k,n}$ & `Sorting map': unilateral $\to$ bilateral &(\ref{eqn:basis_CSI}) \\
\hline
\end{tabular}
\end{center}
\label{tbl:def2}
\end{table}

\subsection{Construction of function spaces and completely orthogonal systems}
\Label{subsec:s41}

In order to introduce the spaces ${\cal H}$ and ${{\cal H}^\Diamond}$, we state two definitions. 

\begin{definition}
\Label{definition:ip}
Define the inner product (among measurable functions on $\mathbb{R} $), parametrized by $k\in\mathbb{Z}$, as 
\begin{eqnarray*}
(f, \,g)_{(k)} := \int_{-\infty}^\infty f(x) \,\overline{g(x)}\ \,(x^2 +1)^k \,dx \, 
\end{eqnarray*}
\end{definition}
\begin{definition}
\Label{definition:sp}
Define the function space 
\begin{eqnarray}
L_{(k)}^2 (\mathbb{R}) :
&=& \Bigl\{ \, f : {\rm measurable }\,\, \Bigl| \,\, \int_{-\infty}^\infty |f(x)|^2 \,(x^2 +1)^k \,dx < \infty \Bigr\}\Label{eqn:sp}
 \\
\Bigl(&=& \bigl\{ \, f : {\rm measurable }\,\, \bigl| \,\, \|f\|_{(k)} < \infty \bigr\} \,\Bigr) . 
\nonumber 
\end{eqnarray}
\end{definition}
Then, 
\begin{eqnarray}
L_{(k)}^2 (\mathbb{R}) \subset L_{(\kappa )}^2 (\mathbb{R} ) \,\,\, {\rm if}\,\,\, k\ge \kappa 
\,\,\,\,\,\,\,\,\,\,\,\,\, {\rm and} \,\,\,\,\,\,\,\,\,\,\, 
L_{(0)}^2 (\mathbb{R}) = L^2 (\mathbb{R}) \, .
\Label{eqn:inclu_sp} \end{eqnarray}

Moreover, obviously, 
\begin{eqnarray}
L_{(k)}^2 (\mathbb{R}) = \, \Bigl\{ \, \frac{f(x)}{(x+i)^k} \,\,\, \Bigl| \,\, f \in L^2(\mathbb{R}) \Bigr\} \,\, . 
\Label{eqn:equiv_sp}\end{eqnarray}
For the spaces ${\cal H}$ and ${{\cal H}^\Diamond}$ introduced for the definition of $\tilde{B}$ in Section \ref{sec:ab}, we will use 
\begin{eqnarray}
{\cal H} & = & L_{(k_0)}^2(\mathbb{R} ) \,\,\, {\rm with } \,\,\, \langle \cdot , \, \cdot \rangle _{{\cal H}} = (\cdot , \, \cdot )_{(k_0)}, 
\Label{eqn:spaces_CSI1}
\\
 {{\cal H}^\Diamond } & = & L_{({k_0^\Diamond })}^2(\mathbb{R} ) \,\,\, {\rm with } \,\,\, \langle \cdot , \, \cdot \rangle _{{{\cal H}^\Diamond }} = (\cdot , \, \cdot )_{({k_0^\Diamond })} , 
\Label{eqn:spaces_CSI2}
\end{eqnarray}
where 
${k_0^\Diamond }\le k_0- s_0$ and $s_0:=\max_{m} \, (\deg p_m -m) $.
Then, Conditions {\bf C1$^+$} and {\bf C2$^+$} trivially hold.

Next, we will introduce the basis function systems $\{e_n\}$ and $\{e_n^\Diamond \}$ for 
these spaces.  
To do this, we need to define the following functions:
\begin{definition}
\Label{definition:ps}
Define the function
\begin{eqnarray}
\psi_{k,\, {\ddot{n}}}(x) := \frac{1}{(x+i)^{k+1}} \left( \frac{x-i}{x+i} \right)^{\ddot{n}} \,\, . 
\Label{eqn:def-psi-CSI}\end{eqnarray}
\end{definition}
Then 
\begin{eqnarray}
\,\,\, \psi_{k,\, {\ddot{n}}}\in L_{(k)}^2(\mathbb{R} ) , \,\, \overline{\psi_{k,\, {\ddot{n}}}(x)} = \psi_{k,\, -{\ddot{n}}-k-1}(x) \,\,\,
 {\rm and } \,\,\, (\psi_{k,\, {\ddot{m}}} \, , \, \psi_{k,\, {\ddot{n}}} )_{(k)} = \pi \,\delta_{{\ddot{m}} {\ddot{n}}} . 
\Label{eqn:properties-psi-CSI}
\end{eqnarray}
The last orthogonal relation is derived easily from calculation of complex integrals by the calculus of residues. 
When $k\ge 0$, as is explained in Section 2 of the paper ~\cite{paper3}, the wavepackets defined by (\ref{eqn:def-psi-CSI}) are `almost-sinusoidally' oscillating wavepackets with spindle-shaped envelopes $|\psi_{k,{\ddot{n}}}(x)| = (x^2 +1)^{-\frac{k+1}{2}} \,$,\, and their approximation (for $\|\cdot\|_{L^2}$) to sinusoidal wavepackets with Gaussian envelopes holds for large $k$. 

For these functions, we have the following lemma, which yields the basis system of our algorithm:
\begin{lemma}
\Label{thm:psi-cons}
$\{\sqrt{\frac{1}{\pi}}\,\psi_{k,\, {\ddot{n}}}\, | \, {\ddot{n}}\in\mathbb{Z} \}$ is an orthonormal basis of $L_{(k)}^2(\mathbb{R} )$.
\end{lemma} \par 
\par The orthonormal property has been shown in the last property of (\ref{eqn:properties-psi-CSI}). Therefore, the proof of completeness in $L_{(k)}^2(\mathbb{R} )$ suffices. This is proved in Appendix \ref{app:pr-cons} from completeness of the Laguerre polynomials, whose details are omitted here, because the Fourier transform of $\psi _{0,{\ddot{n}}}$ can be expressed in terms of the Laguerre polynomial of degree ${\ddot{n}}$. The completeness of $\{\psi _{k,{\ddot{n}}}\,|\, {\ddot{n}}\in\mathbb{Z}\}$ for $k\ne 0$ can therefore also be derived by (\ref{eqn:equiv_sp}) and (\ref{eqn:def-psi-CSI}) ). 

Here we point out some properties of $\psi_{k,\, {\ddot{n}}}$ defined in Definition \ref{definition:sp}, which will be important later.
\begin{theorem}
\Label{thm:ixd-psi}
Any integer $\ddot{n}$ satisfies 
\begin{eqnarray}
\psi_{k,\, {\ddot{n}}}(x) 
&=& - \frac{i}{2} \left( \psi_{k-1,\, {\ddot{n}}} (x) - \psi_{k-1,\, {\ddot{n}}+1} (x), \right)
\Label{eqn:id_CSI} \\
x\, \psi_{k,\, {\ddot{n}}} (x) 
&=& \frac{1}{2} \left( \psi_{k-1,\, {\ddot{n}}} (x) + \psi_{k-1,\, {\ddot{n}}+1} (x), \right)
\Label{eqn:mult_CSI} \\
{\textstyle \frac{d}{dx}}\, \psi_{k ,\, {\ddot{n}}} (x) 
&=& {\ddot{n}} \, \psi_{k+1 ,\, {\ddot{n}}-1} (x) - ({\ddot{n}}+k+1) \, \psi_{k+1 ,\, {\ddot{n}}} (x) .
\Label{eqn:diff_CSI}
\end{eqnarray}
\end{theorem}

This theorem can be derived directly from Definition \ref{definition:ps}. 

These functions are used for the basis systems of ${\cal H}$ and ${{\cal H}^\Diamond }$ as follows: From Lemma \ref{thm:psi-cons}, the following $\{e_n \, \}_{n=0}^{\infty} $ and $\{e_n^\Diamond \}_{n=0}^{\infty} $ are orthonormal basis systems for ${\cal H}$ and ${{\cal H}^\Diamond }$ in 
(\ref{eqn:spaces_CSI1}) and (\ref{eqn:spaces_CSI2}), 
respectively, i.e. Condition {\bf C1} is satisfied:
\begin{eqnarray}
&&e_n= \sqrt{\textstyle\frac{1}{\pi}} \,\psi_{k_0,\, \ddot{n}_{k_0,n}} 
\,\,\, {\rm and} \,\,\,
e_n^\Diamond = \sqrt{\textstyle\frac{1}{\pi}} \,\psi_{{k_0^\Diamond },\, \ddot{n}_{{k_0^\Diamond },n}} 
\Label{eqn:basis_CSI} \\
&&{\rm with} \,\,\, 
{\ddot{n}}_{k,n} := \left\lfloor {\textstyle -\frac{k+1}{2}} \right\rfloor + (-1)^{n+k+1}\, \left\lfloor {\textstyle \frac{n+1}{2}} \right\rfloor,\nonumber
\end{eqnarray}
where $\lfloor a \rfloor$ denotes the largest integer not greater than $a$. 

The indices of functions in $\left\{\,\psi_{k_0,\, \ddot{n}}\, \bigl| \, \ddot{n}\in\mathbb{Z} \right\}$ are bilaterally expressed, while the indices of basis functions in $\{e_n \, \}_{n=0}^{\infty}$ are unilaterally expressed, and they are `matched' to one another by the one-to-one mapping defined by ${\ddot{n}}_{k,n}$ in (\ref{eqn:basis_CSI}). In order to avoid confusion between them, in this paper, the integer indices with double dots\, $\ddot{}$\, denote the bilateral ones in $\mathbb{Z}$, in contrast with the unilateral ones (without double dots) in $\mathbb{Z}^+$. 
For $\ddot{n}$, the order of the above `sorting of the basis' for $\{e_n\}$ is \par ``$-\frac{k_0+1}{2}-\frac{1}{2} ,\, -\frac{k_0+1}{2}+\frac{1}{2},\, -\frac{k_0+1}{2}-\frac{3}{2},-\frac{k_0+1}{2}+\frac{3}{2},\, -\frac{k_0+1}{2}-\frac{5}{2},\, -\frac{k_0+1}{2}+\frac{5}{2},\, ...$'' for even $k_0$, while it is``$-\frac{k_0+1}{2} ,\, -\frac{k_0+1}{2} -1,\, -\frac{k_0+1}{2} +1,\, -\frac{k_0+1}{2} -2,\, -\frac{k_0+1}{2} +2,\, \\ -\frac{k_0+1}{2} -3,\, -\frac{k_0+1}{2} +3,\, ...$'' for odd $k_0$. For $\{e_n^\Diamond \}$, similarly with ${k_0^\Diamond }$ instead of $k_0$. 
The sorting in (\ref{eqn:basis_CSI}) may seem to be somewhat complicated and tricky. However, it is necessary in order to guarantee Conditions {\bf C2} and {\bf C5} later. 
As is mentioned later, other conditions hold,
we can apply the algorithm 
given in Subsection \ref{subsec:s22}
to the quintuplet
$(
P,
L_{(k_0)}^2 (\mathbb{R}),
\{\sqrt{\textstyle\frac{1}{\pi}} \,\psi_{k_0,\, \ddot{n}_{k_0,n}} \}_{n=0}^\infty,
L_{(k_0^{\Diamond})}^2 (\mathbb{R}),
\{\sqrt{\textstyle\frac{1}{\pi}} \,\psi_{{k_0^\Diamond },\, \ddot{n}_{{k_0^\Diamond },n}} \}_{n=0}^\infty
)$.

In addition, when $P(x, \frac{d}{dx})=
\sum_{m=0}^{M}p_m(x)(\frac{d}{dx})^m
$ is a Fuchsian differential operator
with polynomial coefficient functions,
we need to be careful concerning the choice of 
$k_0^{\Diamond}$.
This is because we treat the differential operator $Q(x, \frac{d}{dx})$
instead of $P(x, \frac{d}{dx})$, which is given in (\ref{eq:4-5-1}).
Although 
the differential operator $P(x, \frac{d}{dx})$
requires $k_0^{\Diamond}$ to satisfy
$k_0^{\Diamond}\le k - \max_{m} \, (\deg p_m -m) $,
the differential operator $Q(x, \frac{d}{dx})$
requires $k_0^{\Diamond}$ to satisfy
$k_0^{\Diamond}\le k - \max_{m} \, (\deg p_m -m) -\sum_{n}(M-\mu_n)$,
where $\mu_n$ is the degree of $n$th zero point, whose definition is given in Subsection \ref{subsec:s21}.
With the above condition,
we can apply the algorithm given in Subsection \ref{subsec:s22}
to the quintuplet
$(
Q,
L_{(k_0)}^2 (\mathbb{R}),
\{\sqrt{\textstyle\frac{1}{\pi}} \,\psi_{k_0,\, \ddot{n}_{k_0,n}} \}_{n=0}^\infty,
L_{(k_0^{\Diamond})}^2 (\mathbb{R}),
\{\sqrt{\textstyle\frac{1}{\pi}} \,\psi_{{k_0^\Diamond },\, \ddot{n}_{{k_0^\Diamond },n}} \}_{n=0}^\infty
)$.

\subsection{Check of Conditions {\bf C2}, {\bf C5}, and {\bf C6}}
\Label{subsec:s42}

A recursive use of the relations in Theorem \ref{thm:ixd-psi} results in the following lemma:
\begin{lemma}
\Label{lemma:expansion-recursion}
Let $k_0, j, m \in \mathbb{Z}^+$ and $\kappa \in \mathbb{Z}$. When $\, \kappa \le k_0+m-j $, the function $x^j ({\textstyle \frac{d}{dx}})^m \psi_{k_0,\, {\ddot{n}}}(x)$ can be expressed as a linear combination of $\psi_{\kappa ,\, {\ddot{r}}}(x) $ \par\noindent $({\ddot{r}}={\ddot{n}}-m ,\, {\ddot{n}}-m+1,\, ...\, ,\, {\ddot{n}}+m+k_0-\kappa )\, $ whose coefficients are polynomials of ${\ddot{n}}$ and $k_0$ with degree not greater than $m$. In particular, in this linear combination, the coefficients of the `outermost' terms with $\psi_{\kappa ,\, {\ddot{n}}-m}$ and $\psi_{\kappa ,\, {\ddot{n}}+m+k_0-\kappa }$ are \par\noindent $\, \displaystyle \left(-\frac{i}{2}\right)^{k_0-\kappa -j+m} \left(\frac{1}{2}\right)^j \prod_{t=1}^m ({\ddot{n}}-t+1)$ and $\, \displaystyle \left(\frac{i}{2}\right)^{k_0-\kappa -j+m} \left(\frac{1}{2}\right)^j (-1)^m \prod_{t=1}^m ({\ddot{n}}+t+k_0)$, respectively. 
\end{lemma} \par 
\par The proof is directly derived from Theorem \ref{thm:ixd-psi}, where we apply (\ref{eqn:diff_CSI}) $m$ times, next (\ref{eqn:mult_CSI}) $j$ times and finally (\ref{eqn:id_CSI}) $k_0-\kappa -j+m$ times. Here note that $k_0-\kappa -j+m\ge 0$ from the condition. 
In order to guarantee Conditions {\bf C2} and {\bf C5},
we can derive the following theorem from Lemma \ref{lemma:expansion-recursion}.
Its derivation is based simply on combining the results for the linear combinations of Lemma \ref{lemma:xjDm-banddiag}, which is somewhat complicated and is given in Appendix \ref{app:pr-matrix-banddiag}.

\begin{theorem}
\Label{thm:matrix-banddiag} 
When the coefficient functions $p_m(x)$ $(m=0,1,\ldots ,M)$ are polynomials, the function  
$P(x,{\textstyle \frac{d}{dx}}) e_n(x)$ belongs to ${{\cal H}^\Diamond }$. 
The quantity $b_m^n$ defined in {\bf C2} 
satisfies the following conditions $($a$)$-$($c$)$:

$($a$)$ :\, $ b_m^n=0 \mbox{ if } |m-n|>2M+k_0-{k_0^\Diamond } \,\, .$

$($b$)$ :\, There exists a polynomial $A(x)$ of degree not greater than $M$ such that \par\noindent \hspace{1.5cm} $|b_m^n|\le A(n) $ for any $m,n\in\mathbb{Z}^+$.

$($c$)$ :\, $b_{r-(2M+k_0-{k_0^\Diamond })}^r \ne 0$ for $r\ge 2M+k_0+\max (-{k_0^\Diamond },\, 0)$.
\end{theorem} \par 

This theorem shows that 
the above mentioned quintuplet
satisfies {\bf C2} with $\ell _0= 2M+k_0-{k_0^\Diamond }$.
Hence, 
the dimension $D$ of $V$ is greater than 
the degree $M$ of the differential operator $P$ because $D \ge \ell _0 \ge 2M$.
This theorem also guarantees that 
this quintuplet satisfies {\bf C5} 
with $\ell _0=2M+k_0-{k_0^\Diamond }$ and $j_0=\max ({k_0^\Diamond },\, 0)$ 
if $p_M(\pm i)\ne 0$, 
because $(2M+k_0)-(2M+k_0-{k_0^\Diamond })={k_0^\Diamond }$. 
Hence, the recursive algorithm Theorem \ref{thm:recursive-alg} can be applied when $p_M(\pm i)\ne 0$ and $ ^\forall x\in\mathbb{R} \,\,\, p_M(x)\ne 0$.

The accidental cases where $p_M(i)=0$ or $p_M(-i)=0$ can be easily avoided by a change of coordinate $x \to x+b$ for appropriate $b\in \mathbb{R}$, because $p_M(x)$ has only $M$ roots. More generally, we can use a change of coordinate $x \to ax+b$ for appropriate \par\noindent  $a>0, \, b \in \mathbb{R}$ which is useful not only for this but also for rapid convergence, by `matching' of the scale and the position of the localization between the basis wavepackets and the true solutions. 

Many band-diagonal elements vanish in the matrix $(b_{n}^r)$.
Especially, the equation $b_n^r =0$ holds when $n \le {k_0^\Diamond }-1$ and $r \ge k_0$.
This fact can be shown as follows.
In the above expansion of 
$x^j ({\textstyle \frac{d}{dx}})^m \psi_{k_0,\, {\ddot{n}}}(x)$, 
the terms with 
$\psi_{\kappa ,\, {\ddot{r}}}\,\,$ $\, ({\ddot{r}}\le -1)$ vanish when 
$0\le {\ddot{n}}\le m-1$, and the terms with 
$\psi_{\kappa ,\, {\ddot{r}}}\,\,\, ({\ddot{r}}\ge -\kappa )$ 
vanish when $-k_0-m\le {\ddot{n}} \le -k_0-1$. 
These properties are derived from $ {\textstyle \frac{d}{dx}}\, \psi_{k_0 ,\, 0} (x) = - (k_0+1) \, \psi_{k_0+1 ,\, 0} (x) $ 
(without the term $\ddot{n}\psi_{k_0+1 ,\, \ddot{n}-1}$) 
and ${\textstyle \frac{d}{dx}}\, \psi_{k_0 ,\, -k_0-1} (x) = -(k_0+1) \, \psi_{k_0+1 ,\, -k_0-2} (x) $ (without the term $-(\ddot{n}+k_0+1)\psi_{k_0+1 ,\, \ddot{n}}$) 
which are special cases of (\ref{eqn:diff_CSI}). 
Applying the matching (\ref{eqn:basis_CSI}) to 
the above vanishing property for $\kappa > 0$, 
we obtain the following.
When $n \ge k_0$, the terms in $e_{n'}^\Diamond \,\, (n' \le \kappa -1)$ vanish in this type of expansion of $x^j ({\textstyle \frac{d}{dx}})^m e_n(x)$, which is derived from $\{e_n\, |\, n \le k_0-1 \}=\{\psi_{k_0,{\ddot{n}}}\, |\, -k_0\le {\ddot{n}}\le -1\, \}$ and $\{e_n^\Diamond \, |\, n \le \kappa -1 \}=\{\psi_{\kappa ,{\ddot{n}}}\, |\, -\kappa \le {\ddot{n}}\le -1\, \}$.
So, we conclude that $b_n^r =0$ for $n \le {k_0^\Diamond }-1$ and $r \ge k_0$.

The calculations of $b_{n'}^n := \langle B_P e_n,\, e_{n'}^\Diamond \rangle _{{{\cal H}^\Diamond }}$
need the recursive use of the relations in Theorem \ref{thm:ixd-psi} in the bilateral expression.
Its program can be realized as an integer-type program under {\bf C6} in a practical algorithm explained in the paper ~\cite{paper3}, 
where relations (\ref{eqn:id_CSI}), (\ref{eqn:mult_CSI}) and (\ref{eqn:diff_CSI}) are modularized. 
For the recursion (\ref{eqn:recursive-alg}) in Theorem \ref{thm:recursive-alg}, we have only to know that $b_{n-\ell _0}^{n-2\ell _0 }$, 
$b_{n-\ell _0}^{n-2\ell _0 +1}$, ..., $b_{n-\ell _0}^n$. 
Hence, Condition {\bf C6} holds with this quintuplet.
Here we omit the `sorted version' in the unilateral expression of (\ref{eqn:id_CSI}), (\ref{eqn:mult_CSI}) and (\ref{eqn:diff_CSI}), 
because it is too complicated to use in a practical program.

\subsection{Check of Condition {\bf C3}}
\Label{subsec:s43}
In order to check Condition {\bf C3},
we define 
the operator $\tilde{C}_P$ by 
\begin{eqnarray*}
\left (\tilde{C}_P g \right)(x) 
:= \sum_{m=0}^M \sum_{j=0}^{\deg p_m} \, (-1)^m \overline{p_{m,j}} \, (x^2+1)^{-k_0} ({\textstyle\frac{d}{dx}})^m \Bigl( x^j (x^2+1)^{k_0^\Diamond } \, g(x) \, \Bigr)
\end{eqnarray*}
with $\displaystyle p_m(x):=\sum_{j=0}^{\deg p_m} p_{m,j} \, x^j$ and domain 
\begin{eqnarray*}
D(\tilde{C}_P) =
\{f\in C^M(\mathbb{R}) \cap L_{(k_0^\Diamond )}^2(\mathbb{R} ) 
\,\, |\,\, \tilde{C}_P f \in L_{(k_0)}^2 (\mathbb{R}) \} ,
\end{eqnarray*}
and describe its closed extension by $C_P$.

\begin{theorem}
\Label{thm:basis-in-dom-adj}
Under $k_0^\Diamond \le k_0-s_0$, the operator $C_P$
and 
the operator $B_P$ defined by the action of $P(x,{\textstyle \frac{d}{dx}})$ in Section \ref{sec:ab} satisfy
\begin{eqnarray*}
 ^\forall f \in D(\tilde{B}_P) \,\,\, {\it and } \,\,\, ^\forall n \in \mathbb{Z}\, , \,\,\,\,\, 
\bigl( \, B_P \, f , \, \psi_{k_0^\Diamond ,\, {\ddot{n}}} \, \bigr)_{(k_0-s_0)} 
= \Bigl( \, f, \, C_P \, \psi_{k_0^\Diamond ,\, {\ddot{n}}} \, \Bigr)_{(k_0)} \,\, .
\end{eqnarray*}
\end{theorem} \par 
\par\noindent Theorem \ref{thm:basis-in-dom-adj} guarantees {\bf C3} under the choices 
(\ref{eqn:spaces_CSI1}), 
(\ref{eqn:spaces_CSI2}),
and (\ref{eqn:basis_CSI})
even when $p_M$ has zero points, 
because $D(C_P)$ is dense in $L_{(k_0-s_0)}^2(\mathbb{R} )$.

Since the proof of this theorem requires many pages, it is given in ~\cite{paper2}. Here, we explain briefly the basic idea used in the proof. The equality \par\noindent $\bigl( \, B_P \, f , \, \psi_{k_0^\Diamond ,\, {\ddot{n}}} \, \bigr)_{(k_0-s_0)} 
= \Bigl( \, f, \, C_P \, \psi_{k_0^\Diamond ,\, {\ddot{n}}} \, \Bigr)_{(k_0)} $ 
can be shown by iterative use of the `integration by parts' $\displaystyle \int _a^b p(x)\, q^\prime (x)\, dx = \bigl[ p(x) \, q(x) \bigr]_{x=a}^{x=b} - \int _a^b p^\prime (x)\, q(x)\, dx$ if we can show the disappearance of the contribution of the term $\bigl[ p(x) \, q(x) \bigr]_{x=a}^{x=b} $ at each step of the iteration in the limit as $a\to-\infty$ and $b\to\infty $. We can show its disappearance under the conditions in Theorem \ref{thm:basis-in-dom-adj}, even when $p(x)$ and $q(x)$ do not converge as $x\to\pm\infty $, by means of a `modified kind of smoothing operator' 
$T$ which `blurs' the endpoints $a$ and $b$ so that 
$\bigl(T^np\bigr)(x)$ and $\bigl(T^nq\bigr)(x)$ 
may converge to $0$ as $x\to\pm\infty $ for an integer $n$.

\begin{remark}
\label{rem:nonsym_ex}
The inequality $k_0^\Diamond \le k_0-s_0$ required in Theorem \ref{thm:basis-in-dom-adj} is essential for Condition {\bf C3}. Remember that, even though the ODE can be represented formally by a band-diagonal matrix, this band-diagonal representation is not always valid if {\bf C3} deos not hold. As an example, we consider  the ODE \par\noindent 
$\displaystyle \left(-\frac{i}{2} (x^2 +1)\, \frac{d}{dx}  + (k+1) x -\alpha \right)f(x)=0$
$(k$: integer, $\alpha$: rational constant$)$ with the choice $k_0^\Diamond = k_0 =k$. In this case, $s_0=1$, and hence the inequality $k_0^\Diamond \le k_0-s_0$ is not satisfied. While the solutions $f(x)=C \frac{1}{(x^2+1)^{\frac{k+1}{2}}} \left(\frac{x-i}{x+i}\right)^\alpha $  $(C$: const$)$ belong to  $L^2_{(k)}(\mathbb{R})$, their corresponding vectors do not satisfy the simultaneous linear equations $\displaystyle \sum_{n=\max (0, m-\ell _0 )}^{m+\ell _0 }\!\!\!\!\! b_m^n f_n = 0$ when $\alpha$ is not integer and $C\ne 0$. Hence, due to the contraposition of Theorem \ref{thm:basis-in-dom-adj}, Condition {\bf C3} is not satisfied. 
In other words, in this case, $A_P$ is not symmetric, which has been explained in Remark \ref{rem:nonsym}.
\end{remark}

\subsection{Check of Condition {\bf C4}}
\Label{subsec:s44}

Next, we will show that {\bf C4} is satisfied under the choices 
(\ref{eqn:spaces_CSI1}), (\ref{eqn:spaces_CSI2}), and (\ref{eqn:basis_CSI}). 
When the coefficient functions $p_m$ are polynomial,
the set $S$ of singular points 
is given by the set of zero points of $p_M$, i.e., $p_M^{-1}(0)$.

In a general framework of the theory of elliptic differential equations, 
the following fact is already known.
Assume that a function $f\in L^2(\mathbb{R} )$ satisfies
the conditions $(\frac{d}{dx})^m f\in L^2(\mathbb{R} )$ for $m=1,2...,M-1$ 
and 
belongs to the kernel of the closed extension of an elliptic differential operator 
$\sum_{m=0}^M r_m(x)({\textstyle\frac{d}{dx}})^m$
satisfying $ ^\exists C_1\ge r_M(x) \ge {}^\exists C_2$.
Then, the function $f$ is smooth, i.e., belongs to 
$C^M(\mathbb{R}\setminus p_M^{-1}(0) )\cap L^2(\mathbb{R} )$.
This fact is shown by a generalization to higher order cases of the discussions in~\cite{Gil}, for example.

However, there are many functions that do not satisfy these conditions 
even for true $C^M$-solutions of ODEs. 
For example, the function 
$\displaystyle\,f(x)=\frac{1}{3x^2+1}\, \cos (x^3+x)\,$ is a true solution of the ODE 
\begin{eqnarray*}
\left(\, \bigl(\frac{d}{dx}\bigr)^2 + \frac{6x}{3x^2+1}\, \bigl({\frac{d}{dx}}\bigr) - \frac{6(3x^2-1)}{(3x^2+1)^{2}} - (3x^2+1)^2 \,\right) f(x) = 0 ,
\end{eqnarray*}
which is equivalent with
\begin{eqnarray*}
\left(\, (3x^2+1)^{2} \bigl(\frac{d}{dx}\bigr)^2 
+ 6x(3x^2+1)
\, \bigl({\frac{d}{dx}}\bigr) 
- 6(3x^2-1)(3x^2+1)^{2}- (3x^2+1)^4 \,\right) f(x) = 0 .
\end{eqnarray*}
This solution $f$ can be written as the form
$\sum_{n=0}^\infty f_n \sqrt{\textstyle\frac{1}{\pi}} \,\psi_{k_0,\, \ddot{n}_{k_0,n}} $
with $\{f_n\}_{n=0}^\infty\in  \ell^2(\mathbb{Z}^+)$.

Hence, in order to show {\bf C4},
we should check that
the general solution function $f$ belongs to 
$C^M(\mathbb{R}\setminus p_M^{-1}(0) )\cap L^2(\mathbb{R} )$.
The following theorem is important for check of Condition {\bf C4}.
\begin{theorem}
\Label{thm:nonexistence-pseudo} 
When the coefficient functions
 $p_m(x)$ $(m=0,1,... M)$ are polynomials of $x$
and the linear space defined under 
the quintuplet
$(
P,
L_{(k_0)}^2 (\mathbb{R}),
\{\sqrt{\textstyle\frac{1}{\pi}} \,\psi_{k_0,\, \ddot{n}_{k_0,n}} \}_{n=0}^\infty$,
$L_{(k_0^{\Diamond})}^2 (\mathbb{R}),
\{\sqrt{\textstyle\frac{1}{\pi}} \,\psi_{{k_0^\Diamond },\, \ddot{n}_{{k_0^\Diamond },n}} \}_{n=0}^\infty
)$,
then
for any element $\vec{f}\in  U \cap \ell ^2(\mathbb{Z}^+)$,
there exists
$\varphi \in C^M(\mathbb{R}\backslash p_M^{-1}(0)) $
such that
\begin{eqnarray*} 
\sum_{n=0}^\infty f_ne_n(x_0)
=
\Bigl(\, P(x,{\textstyle\frac{d}{dx}}) \varphi \,\Bigr)
(x_0)=0
\end{eqnarray*}
for $\forall x_0\in \mathbb{R}\backslash p_M^{-1}(0)$.
\end{theorem} 

Theorem \ref{thm:nonexistence-pseudo} seems to provide Condition {\bf C4} directly.
However, this theorem only guarantees the point-wise convergence while 
Condition {\bf C4} requires the convergence in the norm of ${\cal H}$.
Hence, we need the following lemma.

\begin{lemma}
\Label{lemma:conv-norm}
If there exists a function $\varphi \in C^M(\mathbb{R}\setminus p_M^{-1}(0) )$ 
such that \par 
$\displaystyle \lim_{N\to\infty } \sum_{n=0}^N f_ne_n(x)=\varphi (x) $ holds for any $x\in\mathbb{R} $ with a sequence $\{f_n\}_{n=0}^\infty \in \ell ^2(\mathbb{Z}^+)$, then $\displaystyle \lim _{N\to\infty } \Bigl\|\bigl(\sum_{n=0}^N f_ne_n\bigr)-\varphi \, \Bigr\|_{{\cal H}}=0$.
\end{lemma} \par 
\par\noindent\noindent{\em Proof of Lemma }\ref{lemma:conv-norm}: \quad
\par \rm 
Since $\{f_n\}_{n=0}^\infty $ belongs to $\ell ^2(\mathbb{Z}^+)$ and $\{e_n \, \}_{n=0}^{\infty} $ is a CONS of ${\cal H}$, there exists a function $f$ such that $\displaystyle \lim _{N\to\infty } \Bigl\|\bigl(\sum_{n=0}^N f_ne_n\bigr)-f \Bigr\|_{{\cal H}}=0$. Hence, there exists a subsequence $\{N_\nu \}_{\nu =0}^\infty $ such that $\displaystyle \lim_{\nu \to\infty } \sum_{n=0}^{N_\nu } f_ne_n(x)=f(x) \, $ (a.e.). Therefore, from the trigonometric inequality, $\displaystyle \, |f(x)-\varphi (x)| \le \lim_{\nu \to\infty } \Bigl( \, \Bigl| \bigl(\sum_{n=0}^{N_\nu } f_ne_n(x)\bigr)-\varphi (x)\, \Bigr| + \Bigl| \bigl(\sum_{n=0}^{N_\nu } f_ne_n(x)\bigr)-f(x)\, \Bigr| \, \Bigr)=0\, $ (a.e.). Therefore, $\|f-\varphi \|_{{\cal H}}=0$,\, and hence \par\noindent $\displaystyle \lim _{N\to\infty } \Bigl\|\bigl(\sum_{n=0}^N f_ne_n\bigr)-\varphi \, \Bigr\|_{{\cal H}}\le \lim _{N\to\infty } \Bigl(\, \Bigl\|\bigl(\sum_{n=0}^N f_ne_n\bigr)-f \Bigr\|_{{\cal H}} + \|f-\varphi \|_{{\cal H}} \, \Bigr)=
 0$.
\hfill\endproof
\par \par 
\par\noindent Therefore, Lemma \ref{lemma:conv-norm} 
and Theorem \ref{thm:nonexistence-pseudo} guarantee Condition {\bf C4} 
even when $p_M$ has zero points.

\subsection{Application to the non-polynomial case}
\Label{subsec:s41-1}
In this subsection, we explain briefly how the method proposed in this paper can be extended to a more general case where the coefficient functions in the differential operator are not necessarily polynomials but rational functions of $x$.

We can generalize the facts shown in the preceding subsections of this section, 
for differential operators written in the form 
$\displaystyle \, R(x,{\textstyle \frac{d}{dx}})= \sum_{m=0}^M r_m(x) \, ({\textstyle {\textstyle \frac{d}{dx}}})^m$ with rational functions $r_m(x)$ $(m=0,1,...,M)$.
Multiplying the least common multiple of the denominators of $r_m(x)$ ($m=0,1,...,M$),
we obtain a differential operator $P(x,{\textstyle \frac{d}{dx}})
= \sum_{m=0}^M p_m(x) \, ({\textstyle {\textstyle \frac{d}{dx}}})^m$
with the polynomial coefficient functions.
Then, we can apply the algorithm given in Subsection \ref{subsec:s22}
to 
the quintuplet
$(
P,
L_{(k_0)}^2 (\mathbb{R}),
\{\sqrt{\textstyle\frac{1}{\pi}} \,\psi_{k_0,\, \ddot{n}_{k_0,n}} \}_{n=0}^\infty$,
$L_{(k_0^{\Diamond})}^2 (\mathbb{R}),
\{\sqrt{\textstyle\frac{1}{\pi}} \,\psi_{{k_0^\Diamond },\, \ddot{n}_{{k_0^\Diamond },n}} \}_{n=0}^\infty
)$.
Since Condition {\bf C4} holds, 
the numerical result $\tilde{f}$ is close to 
a solution of $P(x,{\textstyle \frac{d}{dx}})f(x)=0$
in $C^M(\mathbb{R} \setminus p_M^{-1}(0)) \cap 
L_{(k_0)}^2 (\mathbb{R})$,
which is 
a solution of $R(x,{\textstyle \frac{d}{dx}})f(x)=0$.

However, any solution of 
$P(x,{\textstyle \frac{d}{dx}})f(x)=0$
in $C^M(\mathbb{R} \setminus p_M^{-1}(0)) \cap 
L_{(k_0)}^2 (\mathbb{R})$
is not necessarily obtained by our algorithm in general.
When the differential operator $R(x,{\textstyle \frac{d}{dx}})$ is Fuchsian,
all of $R(x,{\textstyle \frac{d}{dx}}) f(x)=0$ 
in $
C^M(\mathbb{R} \setminus p_M^{-1}(0)) \cap L_{(k_0)}^2 (\mathbb{R})$
can be approximately obtained by our algorithm.
This fact can be shown as follows.
Any solution of 
$R(x,{\textstyle \frac{d}{dx}}) f(x)=0$ 
in $C^M(\mathbb{R} \setminus p_M^{-1}(0)) \cap L_{(k_0)}^2 (\mathbb{R})$
is a solution of
$Q(x,{\textstyle \frac{d}{dx}}) f(x)=0$,
where the $Q(x,{\textstyle \frac{d}{dx}})$ is given in (\ref{eq:4-5-1}) from 
$P(x,{\textstyle \frac{d}{dx}})$.
Due to Theorem \ref{thm:s35},
any solution of 
$Q(x,{\textstyle \frac{d}{dx}}) f(x)=0$ 
in $C^M(\mathbb{R} \setminus p_M^{-1}(0)) \cap L_{(k_0)}^2 (\mathbb{R})$
can be approximately obtained by our algorithm.
So, we obtain the above fact.

When the ODE $R(x,{\textstyle \frac{d}{dx}})f(x)=0$ has no singular points,
we have the following stronger characterization for the ODE.
This condition is equivalent the non-existence of no zero points in the coefficient function $p_M$.
In this case, all of the solutions of $A_R f=0$ in $L_{(k_0)}^2 (\mathbb{R})$
can be approximately obtained by our algorithm.
This fact can be shown from the following theorem.

\begin{theorem}
Assume that $p_M$ has no zero points.
For any $k\ge 0$, there exists an integer  $k_0^\Diamond$ such that
following conditions for $f \in L_{(k_0)}^2 (\mathbb{R})$ are equivalent.
\begin{description}
\item[(1)]
$A_R f=0$ with 
${\cal H}=L_{(k_0)}^2 (\mathbb{R})$.
\item[(2)]
$B_P f=0$ with 
${\cal H}=L_{(k_0)}^2 (\mathbb{R})$ and 
${\cal H}^{\Diamond}=L_{(k_0^{\Diamond})}^2 (\mathbb{R})$.
\item[(3)]
The $\ell ^2$-sequence 
$\{f_n:=\langle f, e_n\rangle_{{\cal H}}\}_{n=0}^\infty$
belongs to $V$ defined with
with the quintuplet
$(P,
L_{(k_0)}^2 (\mathbb{R}),
\{\sqrt{\textstyle\frac{1}{\pi}} \,\psi_{k_0,\, \ddot{n}_{k_0,n}} \}_{n=0}^\infty,
L_{(k_0^{\Diamond})}^2 (\mathbb{R}),
\{\sqrt{\textstyle\frac{1}{\pi}} \,\psi_{{k_0^\Diamond },\, \ddot{n}_{{k_0^\Diamond },n}} \}_{n=0}^\infty
)$.
\item[(4)]
$P(x,{\textstyle \frac{d}{dx}})f(x)=0$ and  $f \in C^M(\mathbb{R} ) \cap 
L_{(k_0)}^2 (\mathbb{R})$.
\item[(5)]
$R(x,{\textstyle \frac{d}{dx}})f(x)=0$ and  $f \in C^M(\mathbb{R} ) \cap 
L_{(k_0)}^2 (\mathbb{R})$.
\end{description}
\end{theorem}

\begin{proof}
There exist an integer $k_1$ and a constant $c$ such that
$(x^2+1)^{k_1} (\frac{p_M(x)}{r_M(x)})^2 \le c$.
Then, we choose $k_0^\Diamond$
with satisfying the condition 
${k_0^\Diamond }\le \min\{k_0- \max_{m} \, (\deg p_m -m), k_0+ k_1\} $.

The property 
${k_0^\Diamond }\le k_0+ k_1$ yields {\bf (1)}$\Rightarrow${\bf (2)}.
The property 
${k_0^\Diamond }\le k_0- \max_{m} \, (\deg p_m -m)$,
Theorem \ref{thm:c:matrixeq}, and Condition {\bf C4} for the above quintuplet
imply 
{\bf (2)}$\Rightarrow${\bf (3)}$\Rightarrow${\bf (4)}.
Since {\bf (4)}$\Rightarrow$
{\bf (5)}$\Rightarrow$
{\bf (1)} is trivial,
we obtain the desired argument.
\end{proof}

\section{Proofs of Theorems given in Section \ref{sec:ab}}
\Label{sec:CSII}

\subsection{Proof of Theorem \ref{thm:c:matrixeq}}
\Label{subsec:s23}

Now, we prove Theorem \ref{thm:c:matrixeq}, with the following definition, as follows:
\begin{definition}
\Label{definition:sp-tr}
Define 
\begin{eqnarray}
{\cal H}^{(n)}:= {\rm span} (e_0, e_1, ... e_n) 
\,\,\, {\rm and } \,\,\,
{{\cal H}^\Diamond }^{(n)}:= {\rm span} (e_0^\Diamond, e_1^\Diamond, \ldots e_n^\Diamond )
\Label{eqn:def_sptrc}\end{eqnarray}
with 
(\ref{eqn:spaces_CSI1}), (\ref{eqn:spaces_CSI2}), and (\ref{eqn:basis_CSI}), and define the orthogonal projectors $P_n$ and $\tilde{P}_n$ to ${\cal H}^{(n)}$ and ${{\cal H}^\Diamond }^{(n)}$, respectively, with respect to the inner products $\langle \cdot , \, \cdot \rangle_{{\cal H}}$ and $\langle \cdot , \, \cdot \rangle _{{{\cal H}^\Diamond }}$, respectively. 

\end{definition}
\par
{\em Proof of Theorem \ref{thm:c:matrixeq}}
By definition \ref{definition:sp-tr}, $\displaystyle\lim_{n\to\infty} \|P_n f-f\|_{{\cal H}}=0$ for $f\in {\cal H}$, and $\displaystyle\lim_{n\to\infty} \|\tilde{P}_n \tilde{f}-f\|_{{{\cal H}^\Diamond }}=0$ for $\tilde{f}\in {{\cal H}^\Diamond }$. Hence, 
$P_n f$ and $\tilde{P}_n \tilde{f}$
weakly converge to $f$ with respect to the respective inner products. 

The condition {\bf C3} holds for every $\widehat{f}_n$ in any function sequence $\{\widehat{f}_n\in D(\tilde{B})\}$ converging to $f\in D(B)$ 
with respect to ${\cal H}$-norm, and the definition of the graph norm guarantees that $B\widehat{f}_n$ converges to $Bf$ with respect to the ${{\cal H}^\Diamond }$-norm. From these facts, the condition {\bf C3} holds even for 
$f\in D(B) \backslash D(\tilde{B})\, $, and hence it follows that $ ^\forall n, \,\,\, e_n^\Diamond \in D(B^*)$ with 
$B^*=C$. 
Hence, $\langle B(P_mf), \, e_n^\Diamond \rangle _{{{\cal H}^\Diamond }} = \langle P_mf, \, B^*e_n^\Diamond \rangle _{{\cal H}}$, which implies 
\begin{eqnarray*}
\lim_{m\to\infty} \langle B(P_mf), \, e_n^\Diamond \rangle _{{{\cal H}^\Diamond }} = \lim_{m\to\infty} \langle P_mf, \, B^*e_n^\Diamond \rangle _{{\cal H}} = \langle f, \, B^*e_n^\Diamond \rangle _{{\cal H}} = \langle Bf, \, e_n^\Diamond \rangle _{{{\cal H}^\Diamond }} .\end{eqnarray*} 
Therefore, any solution $f\in D(B)$ of $Bf=0$ satisfies 
$\displaystyle \lim_{m\to\infty}\langle B(P_mf), \, e_n^\Diamond \rangle _{{{\cal H}^\Diamond }} = \langle Bf, \, e_n^\Diamond \rangle _{{{\cal H}^\Diamond }} = 0$. 
On the other hand, from {\bf C2}, it is easily shown that \par $ ^\forall m\ge n+\ell _0 , \,\, \tilde{P}_{n}B(P_{m}f)=\tilde{P}_{n}B(P_{n+\ell _0 }f)$. Since $\tilde{P}_{n}\,e_n^\Diamond = e_n^\Diamond $, 
\begin{eqnarray*}
\lim_{m\to\infty}\langle B(P_mf), \, e_n^\Diamond \rangle _{{{\cal H}^\Diamond }} 
&=& \lim_{m\to\infty}\langle B(P_mf), \, \tilde{P}_{n}\,e_n^\Diamond \rangle _{{{\cal H}^\Diamond }}\\ &=& \lim_{m\to\infty}\langle \tilde{P}_{n}B(P_mf), \, e_n^\Diamond \rangle _{{{\cal H}^\Diamond }} 
=
 \langle \tilde{P}_{n}B(P_{n+\ell _0 }f), \, e_n^\Diamond \rangle _{{{\cal H}^\Diamond }} .
\end{eqnarray*}
These facts lead us to $\langle \tilde{P}_{n}B(P_{n+\ell _0 }f), \, e_n^\Diamond \rangle _{{{\cal H}^\Diamond }} =0$, which is equivalent to (\ref{eqn:C1}) in Theorem \ref{thm:c:matrixeq}, 
from {\bf C1} and {\bf C2}, 
because $b_m^n=0$ for $|m-n|>\ell _0$. Thus Theorem \ref{thm:c:matrixeq} holds under {\bf C1}-{\bf C3}.

\hfill\endproof

\vspace{2mm}

\subsection{Proof of Theorems \ref{thm:s35}}
\Label{subsec:s24}

In order to prove Theorem \ref{thm:s35}, we prepare the following definition and lemmata:

\begin{definition}
For a nonnegative integer $L$ and a positive real number $\epsilon$, define the function \begin{eqnarray}\nonumber
W_{L,z,\epsilon }(x):= 1-w_L\left(\frac{|x-z|}{\epsilon }\right) , 
\end{eqnarray} with 
\begin{eqnarray}\nonumber 
\displaystyle w_L(x):=\left\{ \begin{array}{ll} 0 & (\mbox{if }x<0) \\ \\  \displaystyle \left( {\int_0^1 u^L (1-u)^L \, du} \, \right)^{-1} \,  \int_0^x \, u^L (1-u)^L & (\mbox{if }0\le x\le 1) \\ \\ 1 & (\mbox{if } x>1). \end{array}\right. \hspace{5mm} 
\end{eqnarray}

\end{definition}\par\noindent 
This definition results in the following lemma directly. 
\begin{lemma}
\label{lemma:property_W}
The function $W_{L,z,\epsilon }(x)$ has the followimg properties:

$($a$)$:  $W_{L,z,\epsilon }\in C^L (\mathbb{R})$. 

$($b$)$:  $^\forall x \in \mathbb{R}, \quad 0\le W_{L,z,\epsilon} (x)\le 1$ 

$($c$)$:  $^\forall x \in \mathbb{R}\setminus (z-\epsilon, z+\epsilon ), \quad W_{L,z,\epsilon} (x)=1$ 

$($d$)$:  $(\frac{d}{dx})^mW_{L,z,\epsilon }(x)\Bigl|_{x=z}\Bigr.=0$ for $m=0,1,\ldots , L$.

$($e$)$:  $ ^\exists K_{L,m}>0$  such that 
$ ^\forall \epsilon >0$, $^\forall x \in\mathbb{R}$, $ |(\frac{d}{dx})^mW_{L,z,\epsilon }(x)|\le K_{L,m} \epsilon ^{-m}$ 

\hspace{8mm}for $m=0,1,\ldots , L$.
\end{lemma}\par\noindent 
Note that $K_{L,m}$ does not depend on $z$. Lemma \ref{lemma:property_W} and the property of Fuchsian class lead us to the following lemma:
\begin{lemma}
\label{lemma:conv_aboutzero}
Assume that a differential operator $Q(x,\frac{d}{dx})$ 
is Fuchsian, its 
coefficients functions $q_m(x)$ $(m=0,1,\ldots M)$ are 
holomorphic on $\mathbb{R}$,
a Hilbert space ${\cal H}$ satisfies {\bf C1$^+$}.
Choose a solution 
$f\in {C}^M(\mathbb{R}\backslash q_M^{-1}(0))\cap {\cal H}$ 
the ODE  $Q(x,\frac{d}{dx})f=0$ $\bigl($for $x\in\mathbb{R}\backslash q_M^{-1}(0)\bigr)$.
Then, 
for any integer $L$ satisfying $L\ge M$, 
there exist 
positive real numbers $\tilde{\epsilon}$, $\alpha$ and a positive constant $K$ such that 
\begin{eqnarray*}
\displaystyle\, \int _{z_n-\epsilon}^{z_n-0} \left|(x-z_n)^m 
{\textstyle (\frac{d}{dx})}^m \left(W_{L,z_n,\epsilon }(x) f(x) \right)\right|^2 \, dx 
&\le & K \epsilon ^\alpha \\
 \displaystyle\, \int _{z_n+0}^{z_n+\epsilon} \left|(x-z_n)^m {\textstyle (\frac{d}{dx})}^m \left(W_{L,z_n,\epsilon }(x) f(x) \right)\right|^2 \, dx 
&\le & K \epsilon ^\alpha
\end{eqnarray*} 
hold for 
any real number $\epsilon \in (0,\tilde{\epsilon})$,
any zero point $z_n$ of $q_M(x)$, and $m=0,1, \ldots M$.
\end{lemma}\par\noindent 

\par
{\em Proof of Lemma \ref{lemma:conv_aboutzero}}

Corollary \ref{cor:1} implies that 
the set $S$ of the singular points is given by $q_M^{-1} (0)$.
In the following, arrange the zero points so that $z_1<z_2<\ldots z_{N_1}$, and, for a convenience, let $z_0=-\infty$ and $z_{N_1+1}=\infty$ though they are not zero points. 
Since the ODE $Q(x,\frac{d}{dx})f=0$ is of the Fuchsian type, 
as is well known in the theory of the power series expansion about regular singular points~\cite{CoLe}, the solution $f(x)$ can written as $\displaystyle \sum _{\tilde{s}=0}^s (x-z_n)^r \bigl(\log (x-z_n)\bigr)^{\tilde{s}} g_{n,\tilde{s},\pm}(x)$, 
with a holomorphic function $g_{n,\tilde{s},+}(x)$ and $g_{n,\tilde{s},-}(x)$ defined in $(z_n,z_{n+1})$ and $(z_{n-1},z_n)$, 
respectively, about a zero point $z_n$ of $q_M(x)$. 
Here $r$ is the exponent which is a root of the indicial polynomial and $s$ is a non-negative integer not greater than $M-1$.  
Hence, $(\frac{d}{dx})^mf(x)$ can be written as a linear combination  $\displaystyle  \sum _{\tilde{s}=0}^s \sum _{u=0}^m \sum_{v=0}^{\min (m-u,\tilde{s})} c_{u,v} (x-z)^{r-u-v} \bigl(\log (x-z)\bigr)^{\tilde{s}-v} (\frac{d}{dx})^{m-u-v} g_{n,\tilde{s},+}(x)$ 
with coefficients $c_{u,v}$ in $\mathbb{C}$ for $x\in(z_n,z_{n+1})$. 
Moreover, since $g_{n,\tilde{s},+}(x)$ is holomorphic at $x=z_n$ $(n=1,2,\ldots N_1)$, there exist positive real numbers $\tilde{\epsilon }$ and $\tilde{J}_{n,\tilde{m},\tilde{s},+}$ such that 
\begin{eqnarray}\nonumber
^\forall x \in (z_n,z_n+\tilde{\epsilon}), \,\, \left|(\frac{d}{dx})^{\tilde{m}}  g_{n,\tilde{s},+}(x)\right|\le \tilde{J}_{n,\tilde{m},\tilde{s},+}    
\end{eqnarray} 
for $\tilde{m}=0,1,\ldots ,M; \, n=1,2,\ldots N_1;\, \tilde{s}=0,1,\ldots ,s$. 
Similarly, 
there exist positive real numbers $\epsilon _0$ and $\tilde{J}_{n,\tilde{m},\tilde{s},-}$ such that 
\begin{eqnarray}\nonumber
^\forall x \in (z_n-\tilde{\epsilon},z_n), \,\, \left|(\frac{d}{dx})^{\tilde{m}} g_{n,\tilde{s},-}(x) \right|\le \tilde{J}_{n,\tilde{m},\tilde{s},-}    
\end{eqnarray} 
for $\tilde{m}=0,1,\ldots ,M; \, n=1,2,\ldots N_1;\, \tilde{s}=0,1,\ldots ,s$. 

Let $\rho $ be the real part of the exponent $r$ used in the above-mentioned expansion of $f(x)$. 
Since $x(\log x)^{\tilde{s}}$ is infinite-times differentiable for $x\in (0,\infty)$ and $\displaystyle \lim_{x\to +0} x (\log x)^{\tilde{s}}=0$ for $\tilde{s}=0,1,\ldots ,s$, when $t$ is a nonnegative integer, from the above fact, there are positive numbers $\tilde{K}$ and  $\tilde{\epsilon}$ such that 
\begin{eqnarray}\nonumber
^\forall x \in (z_n-\tilde{\epsilon},z_n)\cup (z_n,z_n+\tilde{\epsilon}), 
 ^\forall \tilde{m}\in \{0,1,\ldots m\}, & \, \\ \nonumber \left|{\textstyle (x-z_n)^t(\frac{d}{dx})}^{\tilde{m}} f(x) \right|& \le \tilde{K} |x-z_n|^{t+\rho -\tilde{m}}
\end{eqnarray}
holds for  $n=1,2,\ldots N_1$ if $t$ is a nonnegative integer. 
This fact and the properties (a), (d) and (e) of Lemma \ref{lemma:property_W} lead us to the statement that there are positive numbers $\widehat{K}_{m,\tilde{m}}$ and  $\tilde{\epsilon}$ such that 
\begin{eqnarray}\nonumber
^\forall x \in (z_n-\tilde{\epsilon},z_n)\cup (z_n,z_n+\tilde{\epsilon}),\,\,  ^\forall \epsilon \in (0,\tilde{\epsilon }), && \, \\ 
\nonumber 
\left|
{(x-z_n)^m (\frac{d}{dx})}^m \left(W_{L,z_n,\epsilon }(x) f(x) \right)
\right|
& \le &
 \sum _{\tilde{m}=0}^m \widehat{K}_{m,\tilde{m}} |x-z_n|^{t+\rho -\tilde{m}} \epsilon ^{\tilde{m}-m}
\end{eqnarray} 
holds for $m=0,1,\ldots ,M$ and  $n=1,2,\ldots N_1$. 

Since $\rho >-\frac{1}{2}$ (otherwise, Condition {\bf C1$^+$} would result in $f\notin {\cal H}$), 
the inequality $m>m-\rho -\frac{1}{2}$ holds. Hence, it is easily shown that the lemma holds with $\alpha=\rho +\frac{1}{2}$ $(>0)$, 
because $f$ belongs to ${\cal H}$ and Condition {\bf C1$^+$} is satisfied.  

\hfill\endproof

\par
{\em Proof of Theorem \ref{thm:s35}}

Let $z_1,z_2,\ldots ,z_{N_1}$ $(z_1<z_2<\ldots <z_{N_1})$ be the zero points of $q_M(x)$.  
For a solution $f$ in $C^M(\mathbb{R}\backslash q_M^{-1}(0))$  of the ODE $Q(x,\frac{d}{dx})f=0$, with $\epsilon >0$, define the function 
\begin{eqnarray}\nonumber
f_\epsilon (x):=\left(\prod _{n=1}^{N_1} W_{M,z_n,\epsilon }(x) \right) f(x) .
\end{eqnarray} 
Then, from the properties (b)-(d) of Lemma \ref{lemma:property_W} and Condition {\bf C1$^+$}, it is easily shown that 
\begin{eqnarray}
\label{eqn:conv_H}
\lim _{\epsilon\to +0} \|f_\epsilon (x) -f(x)\|_{\cal H} =0 
\end{eqnarray} 
without any complicated problem caused by the fact $z_n \in q_M^{-1}(0)$. 

On the other hand,  Corollary \ref{cor:1}, 
Lemma \ref{lemma:conv_aboutzero} and Condition {\bf C2$^+$} lead us to 
\begin{eqnarray}
\label{eqn:conv_Hdia}
\lim _{\epsilon\to +0} \left\|\prod _{n=1}^{N_1} Q(x,{\textstyle \frac{d}{dx}})f_\epsilon (x) \right\|_{{\cal H}^\Diamond } =0 . 
\end{eqnarray} 
The convergences (\ref{eqn:conv_H}) and (\ref{eqn:conv_Hdia}) 
imply that $f$ belongs to the domain of $B_Q$ and the equality 
$B_Q f=0$ holds. 
Since $S=q_M^{-1}(0)$,
we obtain the desired argument.
\hfill\endproof
\subsection{Proof of Theorem \ref{thm:conv1}}
\Label{ss43}
The definition of $\ell ^2(\mathbb{Z}^+)$ implies that
\begin{eqnarray}
(\sigma_{K,\infty}^{(\Omega )})^{-1}[0,c 
\underline{\sigma_{K,\infty}^{(\Omega )}} ]
\subset 
V \cap \ell ^2(\mathbb{Z}^+).\Label{eq:1}
\end{eqnarray}
In order to show Theorem \ref{thm:conv1},
we denote the set of normalized vectors in $V$ 
and the $\epsilon $ neighborhood of a normalized vector $\vec{x}$ in this set 
concerning the norm $\|~\|_{\ell ^2,K}$ 
by $O_K$ and $U_{\epsilon,\vec{x}}$.
Thus, from (\ref{eq:1})
it is sufficient to show that
for any $\epsilon >0$,
there exists an integer $N_0$ such that
\begin{eqnarray*}
O_K \cap 
(\sigma_{K,N}^{(\Omega )})^{-1}[0,c 
\underline{\sigma_{K,N}^{(\Omega )}} ]
\subset 
U:=
\bigcup_{\vec{x} \in
O_K \cap 
(\sigma_{K,\infty}^{(\Omega )})^{-1}[0,c 
\underline{\sigma_{K,\infty}^{(\Omega )}} ]} 
U_{\epsilon,\vec{x}}
\end{eqnarray*}
for $N\ge N_0$.

Since $U$ is an open set in $O_K$ and $O_K$ is compact, 
$U^c \cap O_K$ is a compact set in $O_K$.
The relation
\begin{eqnarray*}
\bigcap_{N}
(O_K \cap 
(\sigma_{K,N}^{(\Omega )})^{-1}[0,c 
\underline{\sigma_{K,\infty}^{(\Omega )}} ])
\subset 
O_K \cap 
(\sigma_{K,\infty}^{(\Omega )})^{-1}[0,c 
\underline{\sigma_{K,\infty}^{(\Omega )}} ]
\subset 
U
\end{eqnarray*}
holds.
Taking their complement sets,
we obtain
\begin{eqnarray*}
\bigcup_{N}
(O_K \cap 
(\sigma_{K,N}^{(\Omega )})^{-1}(c 
\underline{\sigma_{K,\infty}^{(\Omega )}} ,\infty))
\supset 
U^c \cap O_K
\end{eqnarray*}
Since
$O_K \cap 
(\sigma_{K,N}^{(\Omega )})^{-1}(c 
\underline{\sigma_{K,\infty}^{(\Omega )}} ,\infty)$ 
is open and 
$O_K \cap 
(\sigma_{K,N}^{(\Omega )})^{-1}(c 
\underline{\sigma_{K,\infty}^{(\Omega )}} ,\infty)
\subset 
O_K \cap 
(\sigma_{K,N+1}^{(\Omega )})^{-1}(c 
\underline{\sigma_{K,\infty}^{(\Omega )}} ,\infty)$,
the compactness of $U^c \cap O_K$
guarantees the existence of an integer $N_0$ such that
\begin{eqnarray*}
O_K \cap 
(\sigma_{K,N}^{(\Omega )})^{-1}(c 
\underline{\sigma_{K,\infty}^{(\Omega )}} ,\infty)
\supset 
U^c \cap O_K
\end{eqnarray*}
for $N \ge N_0$.
Thus, since 
$\underline{\sigma_{K,N}^{(\Omega )}}
\le
\underline{\sigma_{K,\infty}^{(\Omega )}}$,
we obtain
\begin{eqnarray*}
O_K \cap 
(\sigma_{K,N}^{(\Omega )})^{-1}[0,c 
\underline{\sigma_{K,N}^{(\Omega )}}]
\subset 
O_K \cap 
(\sigma_{K,N}^{(\Omega )})^{-1}[0,c 
\underline{\sigma_{K,\infty}^{(\Omega )}}]
\subset U
\end{eqnarray*}
for $N \ge N_0$.

\subsection{Proof of Theorem \ref{thm:conv2}}
\Label{ss44}
In order to show (\ref{eq:2}), it is sufficient to prove that
$(\sigma_{K,\infty}^{(\Omega )})^{-1}[0,
c_0 \underline{\sigma_{K,N}^{(\Omega )}} ]$
contains the subspace
$V \cap \ell ^2(\mathbb{Z}^+)$.

Since $\max_{\vec{x}\in O_K\cap V\cap  \ell ^2(\mathbb{Z}^+)} \|\vec{x}\|_{\ell ^2}$
is finite,
we can choose $c_1$ such that
\begin{eqnarray*}
(\sigma_{K,\infty}^{(\Omega )})^{-1}[0,
c_1 \underline{\sigma_{K,\infty}^{(\Omega )}} ]
=
V \cap \ell ^2(\mathbb{Z}^+).
\end{eqnarray*}
Next, we fix an integer $N_0$ and choose $c_0$ such that
\begin{eqnarray*}
c_0 \underline{\sigma_{K,N_0}^{(\Omega )}} 
\ge
c_1 \underline{\sigma_{K,\infty}^{(\Omega )}} .
\end{eqnarray*}
For any $N \ge N_0$,
\begin{eqnarray*}
c_0 \underline{\sigma_{K,N}^{(\Omega )}} 
\ge
c_0 \underline{\sigma_{K,N_0}^{(\Omega )}} 
\ge
c_1 \underline{\sigma_{K,\infty}^{(\Omega )}} .
\end{eqnarray*}
Thus,
\begin{eqnarray*}
(\sigma_{K,\infty}^{(\Omega )})^{-1}[0,
c_0 \underline{\sigma_{K,N}^{(\Omega )}} ]
\supset
(\sigma_{K,\infty}^{(\Omega )})^{-1}[0,
c_1 \underline{\sigma_{K,\infty}^{(\Omega )}} ].
\end{eqnarray*}
Therefore,
the set 
$(\sigma_{K,\infty}^{(\Omega )})^{-1}[0,
c_0 \underline{\sigma_{K,N}^{(\Omega )}} ]$
contains the subspace
$V \cap \ell ^2(\mathbb{Z}^+)$.

\section{Numerical examples}
\Label{sec:nm}

\begin{figure}[hbt]
\includegraphics{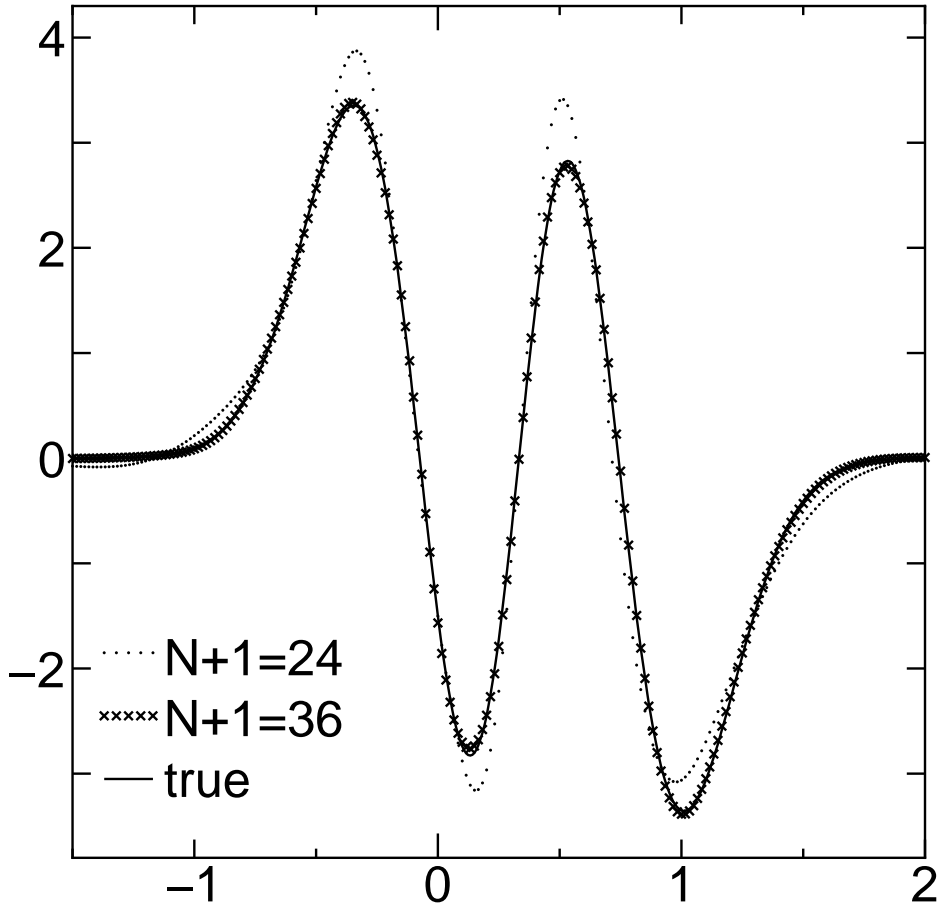}
\caption{Numerical result of functions for the ODE $(\ref{eqn:ex_1}$)}
\label{f2}
\end{figure}

Though the abstract structure of the algorithm is roughly explained in Subsection \ref{subsec:s22} of 
this paper, the detailed explanation of the practical algorithm requires many pages, which is reported in our 
paper~\cite{paper3}, and so we omit it here. Here, we provide some numerical results only. In order that we can observe the accuracy of the algorithm, we chose example ODEs whose exact solutions are known analytically 
and can be written with special functions. 
Though the proposed method can be used widely even for the ODE's which can not be solved analytically, by intent we here give some  examples for ODE's whose true exact solutions can be obtained analytically, in order to clarify how accurate solutions the proposed method gives. 

In the following, we use the bilinear form 
\begin{eqnarray*}
\Omega(\vec{f},\vec{g})=\sum_{n=0}^\infty w_n f_n \overline{g_n}
\end{eqnarray*}
where 
\begin{eqnarray*}
w_n :=\left\{\begin{array}{@{\,}ll} 1 & (n\le K) \\ e^{r(\mu_n-\mu_K)} & \!\!\! (K<n<J) \\ R:=e^{r(\mu_J-\mu_K)} & (n\ge N) \end{array}\right. 
 \,\,\,\mbox{ with } \,\,\, 
\mu_n:=\Bigl|{\ddot{n}}_{k_0,n}-\frac{k_0+1}{2}\Bigr|-\frac{k_0+1}{2} 
\end{eqnarray*}
under the choice $K=2\lfloor \frac{3(N-k_0)}{8} \rfloor + k_0$, 
$J=2\lfloor \frac{7(N-k_0)}{16} \rfloor + k_0$ or 
$K=2\lfloor \frac{7(N-k_0)}{16} \rfloor + k_0$, 
$J=2\lfloor \frac{15(N-k_0)}{32} \rfloor + k_0$, 
and $r=10^8$. The weight number series $\{w_n\}_{n=0}^\infty $ used in this bilinear form may seem to be somewhat complicated, but it is suitable for the symmetry property due to $\overline{\psi_{k_0,{\ddot{n}}}}=\psi_{k_0,-{\ddot{n}}-k-1}$ in (\ref{eqn:properties-psi-CSI}). 

The first example is the third-order ODE 
\begin{eqnarray}
\,\,\,\,\,\,\,\,\,\,\,\,\,\,\, 
f^{\prime\prime\prime} - x f^{\prime\prime}-(81x^2-54x-18\nu )f^{\prime} +\bigl(81x^3-54x^2-(18\nu + 162)x+54\bigr)f=0\, .
\Label{eqn:ex_1}\end{eqnarray}
If $\nu \in\mathbb{Z}^+$, the space of solutions in $C^3(\mathbb{R} ) \cap L_{(k_0)}^2(\mathbb{R} )$ is \par\noindent 
$\{C(\exp \frac{-(3x-1)^2}{2})H_\nu (3x-1) \, |\,C\in \mathbb{C}\} $, where $H_\nu$ is a Hermite polynomial, 
because the differential operator on the left hand side of this ODE can be decomposed as \par\noindent $9\, \bigl({\textstyle\frac{d}{dx}}-x\bigr)\cdot \bigl({\textstyle\frac{1}{3^2}}({\textstyle\frac{d}{dx}})^2-(3x-1)^2 +(2\nu +1)\bigr) 
$
and it can be shown that there is no solution $f$ in $L_{(k_0)}^2(\mathbb{R} )$ such that 
${\textstyle\frac{1}{3^2}}f^{\prime\prime}+\bigl(-(3x-1)^2 +(2\nu +1)\bigr)f$ belongs to $\ker \bigl({\textstyle\frac{d}{dx}}-x\bigr)\backslash \{0\}$. The results with $\nu=3$, $k=4$, $N+1=24, 36$ 
, $K=2\lfloor \frac{3(N-k_0)}{8} \rfloor + k_0$ and $J=2\lfloor \frac{7(N-k_0)}{16} \rfloor + k_0$
are shown in Figure \ref{f2}, 
under  
the normalization $ \, \langle f, \frac{1}{2\pi}(\psi_{k_0,0}+\psi_{k_0,-k_0-1})\rangle_{{\cal H}}= 1$. The errors of the result with $N+1=36$ only are hardly noticeable in this figure.

\begin{table}[h]
\caption{Numerical results for the ratio $\displaystyle \frac{f_2}{f_0}$ under $\nu=0$ and $k=3$ for the ODE $(\ref{eqn:ex_2})$}
\begin{tabular}{|@{}c@{}|@{}c@{}|@{}c@{}|}
\hline 
$N\!\! +\!\! 1$ & ratio $\displaystyle \frac{f_2}{f_0}$ & \, decimal expression of  ratio $\displaystyle \frac{f_2}{f_0}$ \\
\hline 
$\begin{array}{@{\,}ll}  50 \\ \\ \end{array}$ & $ \displaystyle\frac{147826}{391819}  $ & \, $0.3772813467442875409308890074243464456802758\ldots $\\ 
\hline 
$\begin{array}{@{\,}ll}  100 \\ \\ \end{array}$& $\displaystyle\frac{208588565}{552872013}$ & \, $0.3772818303248061138518870912715200145245912\ldots $\\ 
\hline
$\begin{array}{@{\,}ll}  150 \\ \\ \end{array}$& $\displaystyle\frac{1969523740562}{5220298414229}$ & \, $0.3772818303248061138245150519347658988268210\ldots $ \\ 
\hline 
$\begin{array}{@{\,}ll}  200 \\ \\ \end{array}$ & $\displaystyle\frac{531796829098893}{1409547946268876}$ & \, $0.3772818303248061138245150770765762118573286\ldots $\\ 
\hline 
$\begin{array}{@{\,}ll}  250 \\ \\ \end{array}$ & $\displaystyle\frac{651719569462020954}{1727407781341996633}$ & \, $0.3772818303248061138245150770767548665927969\ldots $ \\ 
\hline    
$\begin{array}{@{\,}ll}  300 \\ \\ \end{array}$ & $\displaystyle\frac{150649258697699321707}{399301653535776433703}$ &  \, $0.3772818303248061138245150770767548664028748\ldots $\\ 
\hline \hline
$\begin{array}{@{\,}ll}{\rm true}\!\!\! \\  \\ \end{array}$ & $ \displaystyle 3\! +\! {2\sqrt{2e\pi}\bigl(\!{\rm Erfc}({\textstyle\frac{1}{\sqrt{2}}})\! -\! 1\bigr)} $& \, $0.3772818303248061138245150770767548664028706\ldots $\\ 
\hline 
\end{tabular}
\label{tbl:6.1}
\end{table}
Another example is Weber's differential equation (which is equivalent to 
the Schr\"odinger equation for a harmonic oscillator~\cite{Mes}) 
\begin{eqnarray}
f^{\prime\prime} -x^2 f +(2\nu +1)f=0\, .
\Label{eqn:ex_2}\end{eqnarray}
As is well known, for $\nu\in \mathbb{Z}^+$, the space of solutions in $C^2(\mathbb{R} )\cap L^2(\mathbb{R} )$ is \par\noindent 
$\{C(\exp \frac{-x^2}{2})H_\nu (x) \, |\,C\in \mathbb{C}\} $, which is a subspace of $L_{(k_0)}^2(\mathbb{R} )$ for any $k_0\in\mathbb{Z}^+$. For this example, convergence is very rapid, and we will report its 
accuracy 
by showing within how many digits the ratio between two coefficients $f_n$ and $f_{n'}$ in the expansion $f(x)=\sum_n f_n e_n(x)$ coincides with the true ratio. For example, In Table \ref{tbl:6.1}, we show the results of the ratio $\displaystyle \frac{f_2}{f_0}$ for the case with $\nu=0$, $k_0=3$ 
$K=2\lfloor \frac{7(N-k_0)}{16} \rfloor + k_0$ and $J=2\lfloor \frac{15(N-k_0)}{32} \rfloor + k_0$, 
where the true ratio is obtained 
analytically (not numerically)  
by means of the computer algebra software package ``Methematica". Similar accuracy is observed for other ratios between the coefficients with small $n$ and $n'$. With $N+1=7000$, we obtained a result where it coincided with the true value up to 
$340$ digits. 
In Figure \ref{fig:ho}, we plot how the number of 
significant digits of this ratio depends on $N$. Moreover, we found that the rational ratios obtained in this case have almost a `full precision', because the proportion  
\begin{eqnarray}
\hspace{5mm} \rho :=\frac{(\mbox{number of significant digits of the ratio})}{(\mbox{number of digits of numerator})+(\mbox{number of digits of denominator})} 
\Label{eqn:def_ratio_digits}\end{eqnarray}
almost equals $1$ for $N+1\ge 100$ as is shown in Figure \ref{fig:rho}. (In this case, the ratio $\frac{f_2}{f_0}$ has no imaginary part due to a symmetry.)

\begin{figure}[t]
\begin{minipage}{0.5\linewidth}
\includegraphics{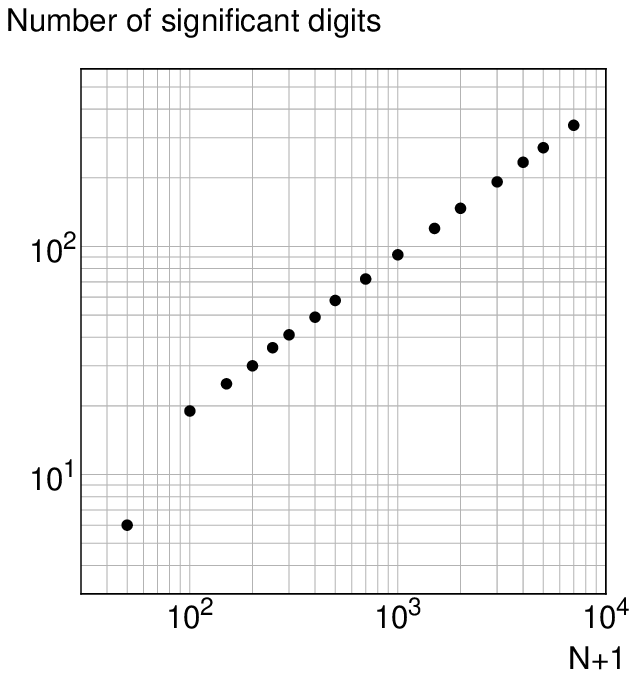}
\caption{Number of significant digits of the ratio $\frac{f_2}{f_0}$, for ODE $(\ref{eqn:ex_2})$}
\Label{fig:ho}
\end{minipage}
\begin{minipage}{0.48\linewidth}
\mbox{Proportion $\rho$ defined in (\ref{eqn:def_ratio_digits})}
\includegraphics{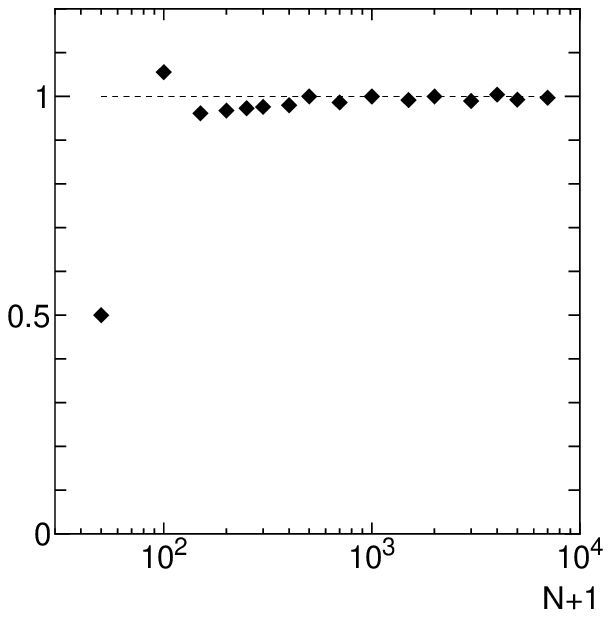}
\caption{Proportion $\rho$ defined in $(\ref{eqn:def_ratio_digits})$, for ODE $(\ref{eqn:ex_2})$}
\Label{fig:rho}
\end{minipage}
\end{figure}

\begin{figure}[t]
\begin{minipage}{0.51\linewidth}
\includegraphics{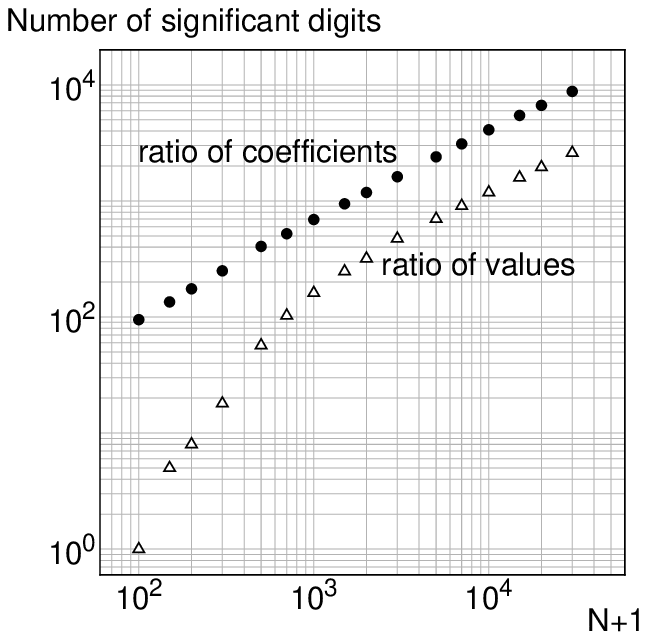}
\caption{Number of significant digits of the ratio $\frac{f_7}{f_5}$ and Number of significant digits of the ratio $\frac{f(1/30)}{f(0)}$, for ODE $(\ref{eqn:ex_3})$}
\Label{fig:ho30}
\end{minipage}
\begin{minipage}{0.48\linewidth}
$\mbox{Proportion $\rho$ defined in (\ref{eqn:def_ratio_digits})}$ 
\includegraphics{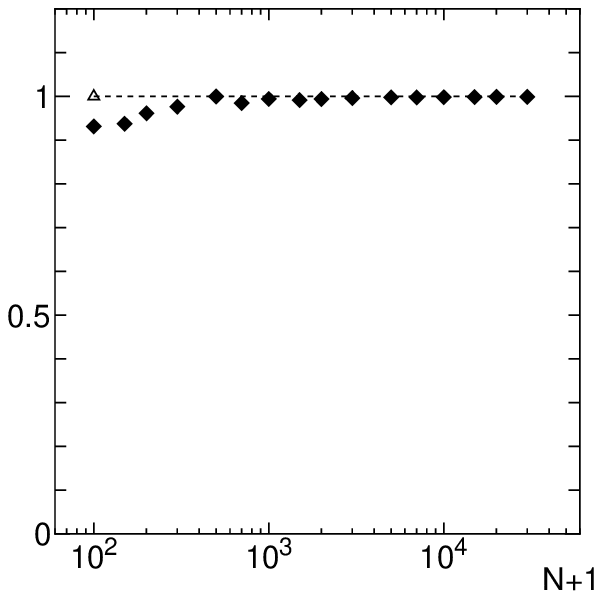}
\caption{Proportion $\rho$ defined in $(\ref{eqn:def_ratio_digits})$, for ODE $(\ref{eqn:ex_3})$}
\Label{fig:rho30}
\end{minipage}
\end{figure}

Moreover, under the scale change $x\to 30x$, the accuracy is 
improved 
very much. Under this scale change, ODE (\ref{eqn:ex_2}) is modified to 
\begin{eqnarray}
\frac{1}{(30)^2}f^{\prime\prime} -(30)^2x^2 f +(2\nu +1)f=0\, .
\Label{eqn:ex_3}\end{eqnarray}
The results for this ODE with 
$\nu=0$, $k_0=6$, $K=2\lfloor \frac{7(N-k_0)}{16} \rfloor + k_0$ and $J=2\lfloor \frac{15(N-k_0)}{32} \rfloor + k_0$
are given in Fig \ref{fig:ho30}. In this case, with $N\ge 100$, the number of the significant digits between two coefficients is infinite (i.e. perfectly exact ratio is obtained) when the true ratio is rational which occurs for the ratio among $f_0, f_1, \ldots , f_5$, and it is very large even when the true ratio is irrational. 
For example, for the ratio $\frac{f_7}{f_5}$ (which is irrational), the ratio obtained numerically by the proposed method coincides within 8783 digits to the true ratio when $N+1=30000$.  
Moreover, for the ratio between the values of the solution function at two points $\frac{f(1/30)}{f(0)}$, there the numerical result by the proposed method coincides within 2599 digits to the true ratio. (The number of significant digits seems to be proportional to $N^{0.79}$ empirically in the case when $N$ is sufficiently large.)  As for the ratio defined in (\ref{eqn:def_ratio_digits}), the numerical results by the proposed method give almost a `full precision' 
al so in this case, which is shown in Fig \ref{fig:rho30}. There results show evidently how accurate the proposed method is. 

Next, we give examples in Fuchsian class where $p_M(x)$ has zero points. 
Here, we show how numerical results converge to true solutions in 
such a case where the true solutions are written with the associate Legendre 
function. In this case, 
we are successful to extract only the true solutions defined only in 
$(-1,1)$ 
by the proposed method, where the obtained solutions is almost zero outside these intervals. 

\begin{figure}[t]
\includegraphics{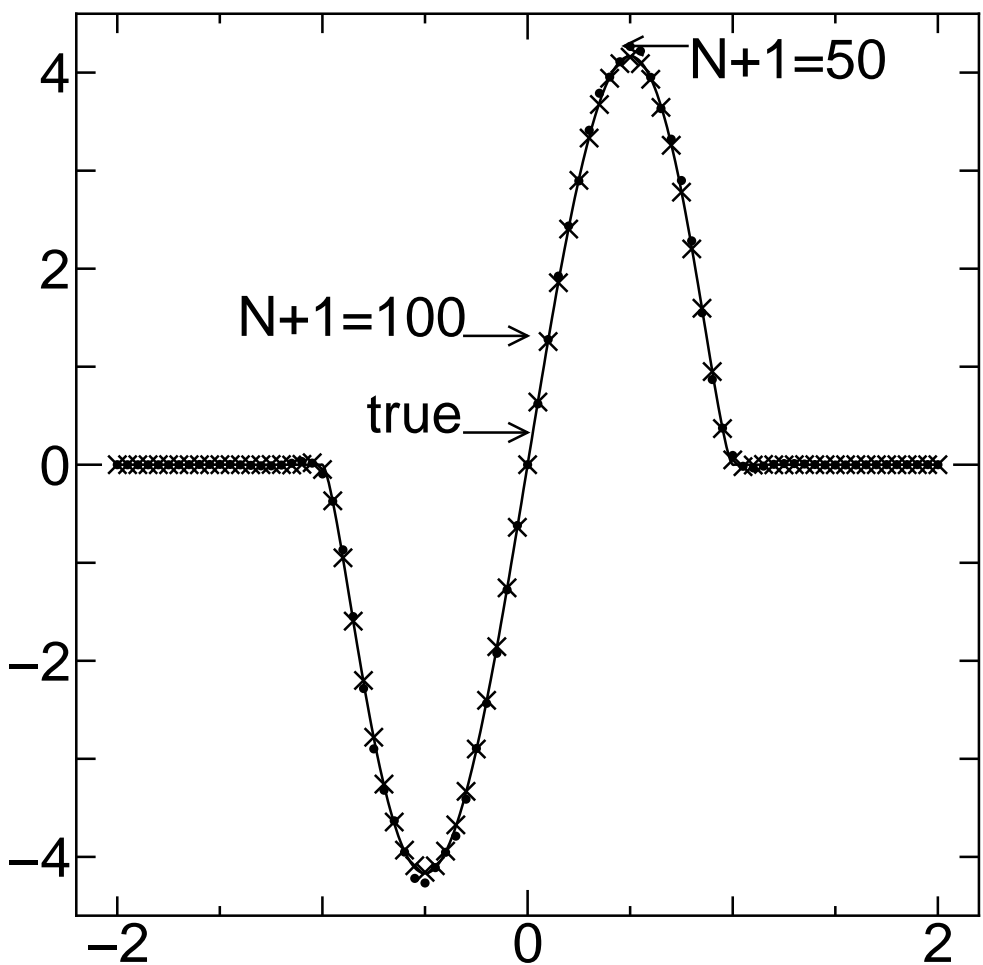}
\caption{Numerical result of functions for the associate Legendre differential equation  $(\ref{eqn:ex_4}$)}
\Label{fig:assleg}
\end{figure}

An example of such a case is for the associate Legendre differential equation
\begin{eqnarray}
(1-x^2)f^{\prime\prime}-2xf^\prime+\left(\nu(\nu +1)-\frac{\mu ^2}{1-x^2}\right)f=0.
\Label{eqn:ex_4}\end{eqnarray}
By means of the discussion in Subsection \ref{subsec:s41-1}, we can treat this ODE by the proposed algorithm as the 
Fuchsian-type 
ODE 
\begin{eqnarray}
(1-x^2)^2f^{\prime\prime}-2x(1-x^2)f^\prime+\left(\nu(\nu +1)(1-x^2)-\mu ^2\right)f=0 
\Label{eqn:ex_5}\end{eqnarray}
whose coefficient functions are polynomials. As is well known, there are three intervals $(-\infty ,-1)$, $(-1,1)$ and $(1,\infty )$ within which smooth solutions are defined, because the coefficient function $p_2(x)$ of the highest order term has two zero points $x=\pm 1$. However, none of the solutions defined in the intervals $(-\infty ,-1)$ and $(1,\infty )$ is  square-integrable, and hence the space of solutions in $L_{(k_0)}^2(\mathbb{R})$ ($\subset L^2(\mathbb{R})$) is the one-dimensional space $\{C \cdot 1_{[-1,1]}(x)\cdot  (1-x^2)^\frac{\mu}{2} L_\nu^\mu (x) \, | \, C\in \mathbb{C} \}$\,\, ($1_I(x)$: indicator function, $(1-x^2)^\frac{\mu}{2} L_\nu^\mu (x)$: associate Legendre function). The results for this ODE with 
$\mu=3$, $\nu=4$, $k_0=6$, $K=2\lfloor \frac{7(N-k_0)}{16} \rfloor + k_0$, $J=2\lfloor \frac{15(N-k_0)}{32} \rfloor + k_0$ and $N+1=50,100$ 
are given in Fig \ref{fig:assleg}. Note that there the solutions are normalized by $ \, \langle f, \frac{1}{2\pi}(\psi_{k_0,0}+\psi_{k_0,-k_0-1})\rangle_{{\cal H}}= 1$. Surprisingly, almost only the component in $\{C \cdot 1_{[-1,1]}(x)\cdot  (1-x^2)^\frac{\mu}{2} L_\nu^\mu (x) \, | \, C\in \mathbb{C} \}$ is `automatically' extracted, and the numerical solutions are almost zero outside the interval $(-1, 1)$, in spite of the existence of singularities at $x=\pm 1$. However, the convergence to the true solution is not so rapid as the cases where $p_M(x)$ has no zero points, though it converges to the true solution anyway.

For the ODEs whose exact solutions can be written by the (associate) Laguerre functions $x^{\frac{\mu}{2}}  e^{-\frac{x}{2}} L_\nu^\mu(x)$ within the interval $(0,\infty )$, we have already had similar results to this, where the obtained numerical solutions are almost zero for $x<0$.

\section{Discussion}
\Label{sec:di}

\subsection{Some properties of the basis functions used in this study}

The basis systems $\{e_n \, \}_{n=0}^{\infty} $ and $\{e_n^\Diamond \}_{n=0}^{\infty} $ are closely related to Fourier series, by the change of variable $\theta = 2 \arctan x$, as is shown in subsection {\bf 2.4} of the paper ~\cite{paper2}. 
(The same change of variable has been used for a description of analytic unit quadrature signals with nonlinear phase~\cite{Qia}~\cite{Che}.) 

The function $\psi_{k,\, 0}$ is identical to the Cauchy wavelet~\cite{Hol} used for continuous wavelet transformation~\cite{Dau}. Moreover, when $k$ is even, $\psi_{k,\, {\ddot{n}}}$ is closely related to the number state associated with $\mathfrak{su}(1,1)$ in a representation of $\mathfrak{su}(1,1)$ which can be formulated by adding a third generator to the two generators of the $ax+b$ group~\cite{SaHa}.

\subsection{Extension to inhomogeneous differential equations}
The algorithm proposed in this paper is easily extended to linear inhomogeneous ordinary differential equations with inhomogeneous terms in ${{\cal H}^\Diamond }$. This extension only requires substitution of the right hand side $0$ of the simultaneous linear equations $\sum_n b_m^n f_n =0$ $(m\in Z^+ )$ by the $\widetilde{H}$-inner-products between the inhomogeneous term
 and the basis function $e_m^\Diamond $.

\subsection{Modification of the method for the eigenvalue-eigenvector problem}
We have already proved that the proposed method can be applied for  eigenfunction problems of self-adjoint operators with given eigenvalues, under some conditions, which will be reported in another paper~\cite{paper4}. 
In order to apply the proposed method to the eigenvalue-eigenvector problem for a linear 
operator, we must have a method to obtain the eigenvalues, because the eigenvalue is regarded as a fixed parameter of the characteristic equation in the proposed method. In the case of discrete eigenvalues, if an eigenvalue is not exact, the function satisfying the characteristic equation does not belong to ${\cal H}$, and hence its corresponding vector is not square-summable. 

However, when we truncate the algorithm within a finite number of dimensions, the square-summability is not distinguishable. The number sequence obtained by our method for an approximate eigenvalue decays within a finite number of dimensions as rapidly as the number sequence corresponding to the true eigenvector. As the approximation of the eigenvalue is better, it decays for more dimensions. From this fact, we can propose a method to find the eigenvalue by observing the location of the bottom of the valley of the ratio $\sigma_{K}^{(\Omega)}$. 
Here we give an example of such 
valleys in Figure \ref{fig:val}.
In this example, we are successful to separate two eigenvalues which are very contiguous by the `tunnel effect', for a Schr\"{o}dinger equation with quantum-double-well-type potential function.

\begin{figure}[h]
\begin{minipage}{0.48\linewidth}
\begin{picture}(100,155)
\put(0,0){\scalebox{.5}{\includegraphics{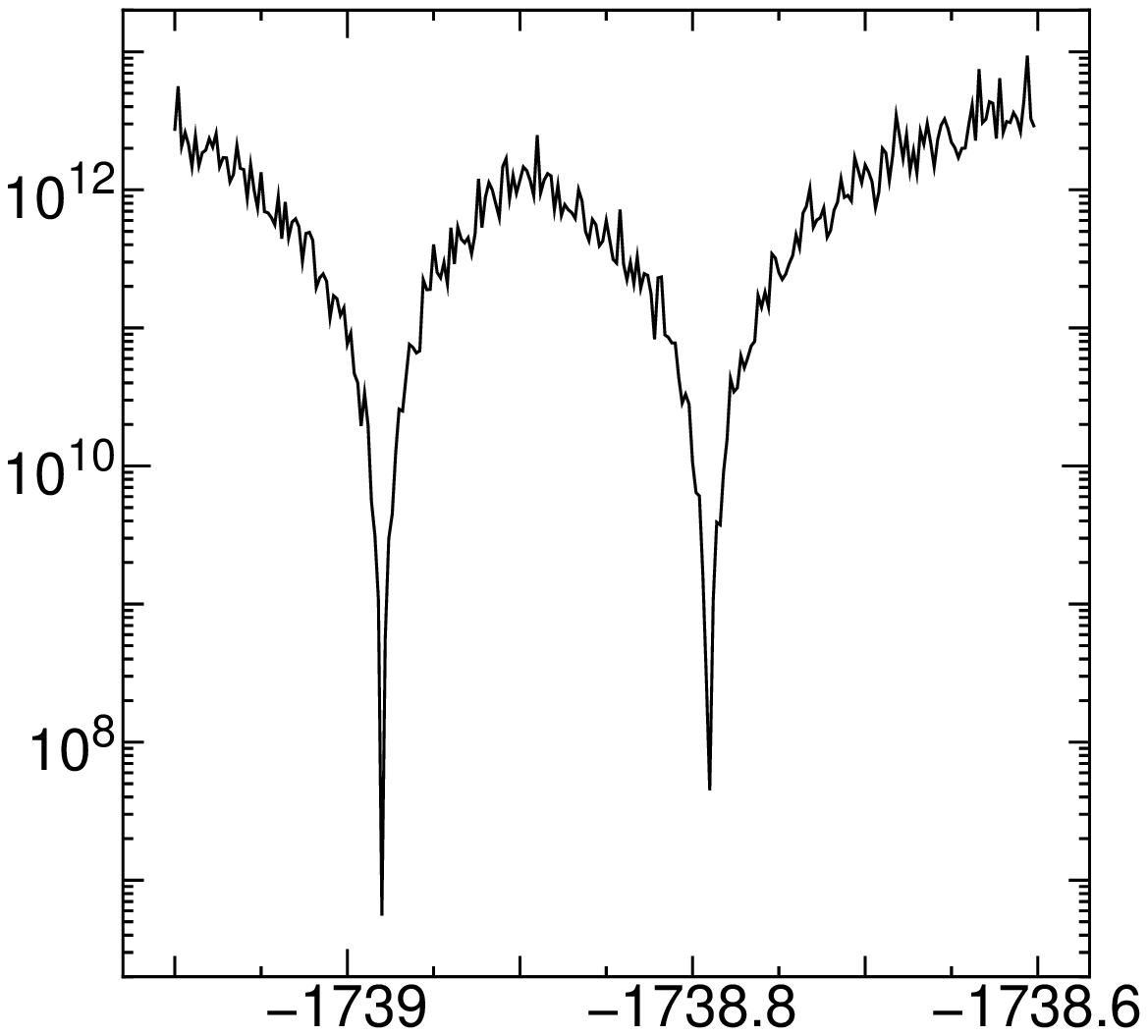}}}
\end{picture}
\caption{Example of `valleys' of the ratio between the two norms}
\Label{fig:val}
\end{minipage}
\begin{minipage}{0.48\linewidth}
\begin{picture}(100,155)
\put(0,0){\scalebox{.35}{\includegraphics{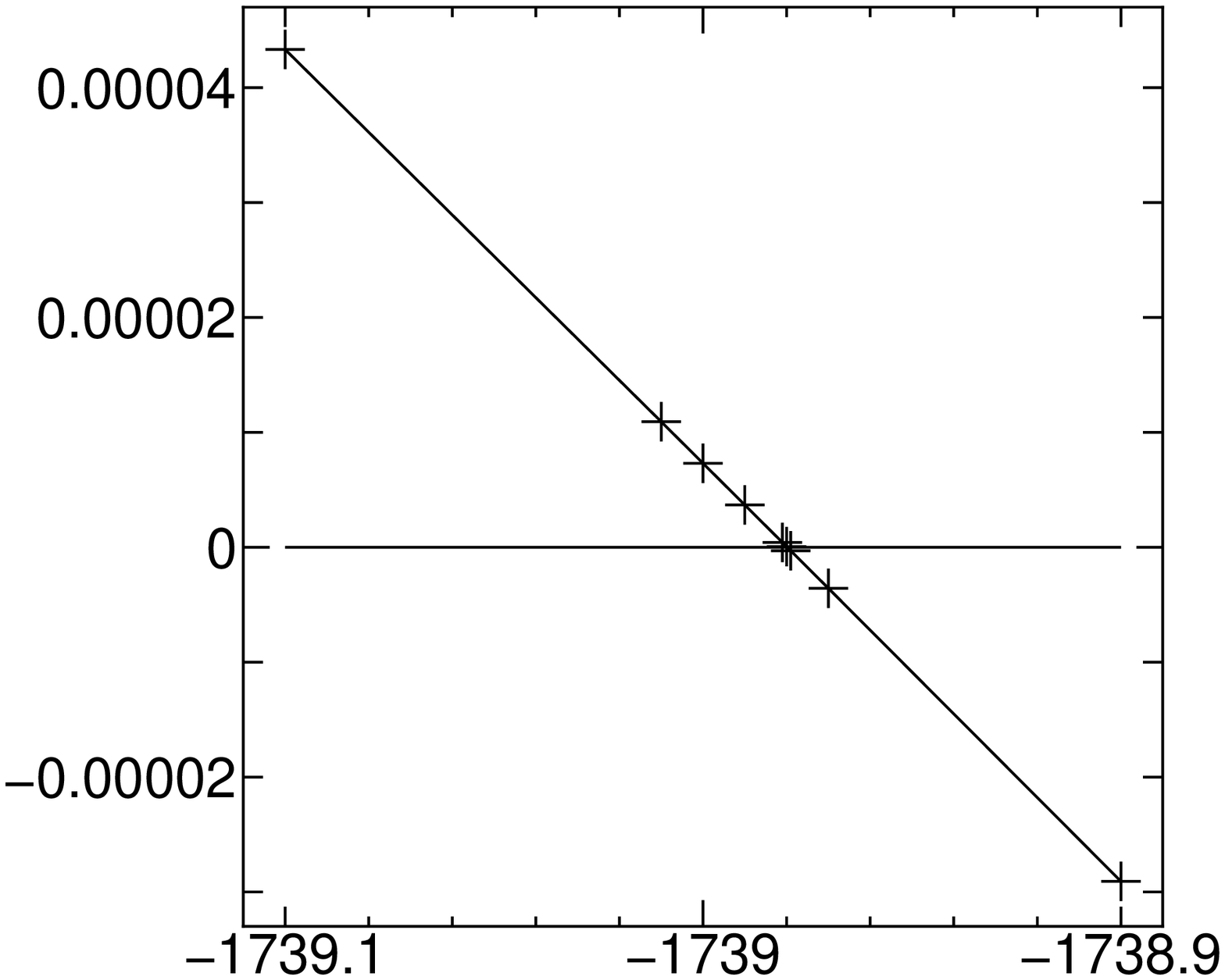}}}
\end{picture}
\caption{Possibility of a more precise interpolation by means of an index almost linear to the deviation}
\Label{fig:lin}
\end{minipage}
\end{figure}


Moreover, we have already invented another faster 
and more effective 
method for finding eigenvalues in a very high accuracy, based on a more analytical idea. 
This idea utilizes linear interpolations by means of some indices almost linear to the deviation of the eigenvalue 
(see Fig.\ref{fig:lin}) which are calculated directly from numerical results. 

\subsection{Possibility of the extension to partial differential equations}
A similar idea to the proposed method can be applied to linear partial differential equations. 
However, the number of linearly independent solutions of simultaneous linear equations is not fixed but increasing as $N$ increases for linear partial differential equations, while it is fixed at $p=j_0+\ell _0 -1$ for linear ordinary differential equations. Therefore, we have to estimate how much memory and how many calculations would be required.

\subsection{Possibility of the extension to weakly non-linear differential equations}
This algorithm has the possibility of extension to nonlinear differential equations because of the following properties:
>From the definition of $\psi_{k_0, {\ddot{n}}}$, \, the relation $\psi_{k_0,{\ddot{n}}_1}(x) \cdot \psi_{k_0,{\ddot{n}}_2}(x) = \psi_{2k_0+1,\, {\ddot{n}}_1+{\ddot{n}}_2}(x)$ holds. The combination of this fact and Lemma \ref{lemma:expansion-recursion} results in the fact that the product $\psi_{k_0,{\ddot{n}}_1}(x) \cdot \psi_{k_0,{\ddot{n}}_2}(x)$ can be expressed as a linear combination of $\psi_{k_0,\, {\ddot{n}}_1+{\ddot{n}}_2}(x), \psi_{k_0,\, {\ddot{n}}_1+{\ddot{n}}_2+1}(x), \ldots , \psi_{k_0,\, {\ddot{n}}_1+{\ddot{n}}_2+k_0+1}(x)$. Similarly, the product of more than three basis functions can be written as a linear combination of finite numbers of the same basis functions. If the nonlinearity is weak, we can apply the proposed method to the successive approximation method for nonlinear differential equations, because of this property. However, for the nonlinear case, it is more difficult to find a proof of convergence and an upper bound for errors, than it is for the
  linear case.

\section{Conclusions}
\Label{sec:cc}

We have proposed an integer-type algorithm which can determine 
accurately a basis system 
for the space of solutions in 
${\cal H}$ 
of the $M$-th order ODE $\displaystyle\Bigl( \sum_{m=0}^M p_m (x) \bigl({\textstyle \frac{d}{dx}}\bigr)^m \Bigr) f(x) = 0$ with polynomials or rational functions for the coefficient functions $p_m$ $(m=0,\ldots ,M)$ under certain conditions. The basic structure of this algorithm has been shown in a more general framework and several conditions have been stated for the validity of this structure. Next, we have provided choices for the spaces and their basis systems satisfying these conditions, with detailed checks of these conditions. Thus, the validity of the proposed method has been proved. 

Moreover, we have shown convergence of the results of this method to true solutions of the differential equations, under the conditions required for the structure of the algorithm.
Numerical results have indicated that this method has high accuracy. We have provided examples to show how the results converge to true solutions as the dimension of the subspace increases. 

This method will be extended or generalized for inhomogeneous equations, partial equations and weakly nonlinear equations in the near future, as has been mentioned in Section \ref{sec:di}. Analyses of the accuracy and the 
amount  
of calculations required are also future problems. Moreover, it is our intent to apply this method, with some modifications, to the scattering problem in quantum mechanics.

\section*{Acknowledgments}
MH was partially supported by MEXT through a Grant-in-Aid for Scientific Research in the Priority Area "Deepening and Expansion of Statistical Mechanical Informatics (DEX-SMI)", No. 18079014 and
a MEXT Grant-in-Aid for Young Scientists (A) No. 20686026. 
The Center for Quantum Technologies is funded by the Singapore Ministry of Education
and the National Research Foundation as part of the Research Centres of Excellence
programme.

\appendix

\section{Proof of Lemma \ref{thm:psi-cons}}
\Label{app:pr-cons}
\par\noindent\noindent{\em Proof of Lemma }\ref{thm:psi-cons}: \quad 
 \rm 
>From the last property of (\ref{eqn:properties-psi-CSI}), $\{\sqrt{\frac{1}{\pi}}\,\psi_{k,\, {\ddot{n}}}\, | \, {\ddot{n}}\in\mathbb{Z} \}$ is orthonormal. Therefore, we have only to prove  the completeness in $L_{(k)}^2(\mathbb{R} )$. Let ${\cal F}$ be the Fourier transformation, where the Fourier transform of a function $f$ is denoted by $\displaystyle \bigl({\cal F}f\bigr)(y):=\frac{1}{\sqrt{2\pi }}\int_{-\infty}^\infty f(x) e^{-iyx} \, dx$. Some calculations by residue calculus result in 
\begin{eqnarray}
 ^\forall 
\ddot{n}
\ge 0, \,\,\,\ ({\cal F}\psi_{0,\, {\ddot{n}}} )(y) = \left\{ \begin{array}{@{\,}ll} \displaystyle i \, \sqrt{2\pi} \,\,\, e^{-y} \, L_{\ddot{n}}(2y) & (y\ge 0) \\ \\ \displaystyle 0 & (y<0) \end{array} \right.
\Label{eqn:freq0}\end{eqnarray}
where $L_n(x)$ denotes the Laguerre polynomial of degree $n$. On the other hand, since $\overline{\psi_{0,\, {\ddot{n}}}(x)}= \psi_{0,\, -{\ddot{n}}-1}(x)$ from (\ref{eqn:properties-psi-CSI}), a property of the Fourier transform leads us to 
\begin{eqnarray}
 ^\forall 
\ddot{n}
\ge 0, \,\,\,\ ({\cal F}\psi_{0,\, -{\ddot{n}}-1} )(y) = \left\{ \begin{array}{@{\,}ll} \displaystyle -i \, \sqrt{2\pi} \,\,\, e^{\, y} \, L_{\ddot{n}}(-2y) & (y\le 0) \\ \\ \displaystyle 0 & (y>0). \end{array} \right. 
\Label{eqn:freq0-}\end{eqnarray}
Here, let 
\begin{eqnarray*} 
{\cal L}_{(0)}^- &:=& \bigl\{\, \sum_{{\ddot{n}}=-\infty}^{-1} \xi_{\ddot{n}} \,\psi_{0,\, {\ddot{n}}}(x) \,\,\bigl| \,\, \xi_{\ddot{n}} \in {\mathbb{C}}, \,\, \{ \xi_{\ddot{n}} \} \in \ell ^2(\mathbb{Z}\backslash \mathbb{Z}^+ ) \,\bigr\},
\\ 
{\cal L}_{(0)}^+ &:=& \bigl\{\, \sum_{{\ddot{n}}=0}^\infty \xi_{\ddot{n}} \,\psi_{0,\, {\ddot{n}}} (x) \,\,\bigl| \,\, \xi_{\ddot{n}} \in {\mathbb{C}} , \,\, \{ \xi_{\ddot{n}} \} \in \ell ^2(\mathbb{Z}^+ ) \,\bigr\} \,\, .
\end{eqnarray*}
Then, from the well-known fact that the set 
$
\{\, e^{-\frac{t}{2}} \, L_n(t), \,\, t\ge 0 \,\, |\,\, n \in \mathbb{Z}^+ \} 
$
is complete in 
$
L^2 \bigr(\, \mathbb{R}^+ \,\bigr) \, , 
$
we can show that $\{ \, {\cal F}f \,\, |\, f\in {\cal L}_{(0)}^+ \} = L^2 \bigr(\, \mathbb{R}^+ \,\bigr)$. Similarly, from (\ref{eqn:freq0-}) and this fact, $\{ \, {\cal F}f \,\, |\, f\in {\cal L}_{(0)}^- \} = L^2 \bigr(\, \mathbb{R}^- \,\bigr)$. Since the null functions in $L^2(\mathbb{R} )$ which are nonzero only at $y=0$ in the frequency domain belong to the kernel of the inverse Fourier transformation, from the Planchrel theorem, 
\begin{eqnarray}
L^2(\mathbb{R}) = {\cal L}_0^- \oplus {\cal L}_0^+ 
\, =\, \bigl\{\, \sum_{{\ddot{n}}=-\infty}^{\infty } \xi_{\ddot{n}} \,\psi_{0,\, {\ddot{n}}}(x) \,\,\bigl| \,\, \xi_{\ddot{n}} \in {\mathbb{C}}, \,\, \{ \xi_{\ddot{n}} \} \in \ell ^2(\mathbb{Z} ) \,\bigr\} \,\, , 
\Label{eqn:equiv_sp_CSI}
\end{eqnarray}
and hence $\{\sqrt{\frac{1}{\pi}}\, \psi_{0,{\ddot{n}}}\, |\, {\ddot{n}}\in \mathbb{Z} \}$ is complete in $L^2(\mathbb{R})=L_{(0)}^2(\mathbb{R})$.
Then, since $\displaystyle \,\psi_{k,\, {\ddot{n}}}(x) = \frac{\psi_{0,\,{\ddot{n}}}(x)}{(x+i)^k} \, $, from (\ref{eqn:equiv_sp}) and (\ref{eqn:equiv_sp_CSI}), 
\begin{eqnarray*} 
L_{(k)}^2(\mathbb{R}) 
&=& \bigl\{\, \frac{1}{(x+i)^k} \sum_{n=-\infty}^{\infty } \xi_n \,\psi_{0,\, {\ddot{n}}}(x) \,\,\bigl| \,\, \xi_{\ddot{n}} \in {\mathbb{C}}, \,\, \{ \xi_{\ddot{n}} \} \in \ell ^2(\mathbb{Z} ) \,\bigr\} \\ 
&=& \bigl\{\, \sum_{{\ddot{n}}=-\infty}^{\infty } \xi_{\ddot{n}} \,\psi_{k,\, {\ddot{n}}}(x) \,\,\bigl| \,\, \xi_{\ddot{n}} \in {\mathbb{C}}, \,\, \{ \xi_{\ddot{n}} \} \in \ell ^2(\mathbb{Z} ) \,\bigr\} \,\, , 
\end{eqnarray*}
and hence $\{\sqrt{\frac{1}{\pi}}\, \psi_{k,{\ddot{n}}}\, |\, {\ddot{n}}\in \mathbb{Z} \}$ is complete in $L_{(k)}^2(\mathbb{R})$. 
\hfill\endproof
\par

\section{Proof of Theorem \ref{thm:matrix-banddiag}}
\Label{app:pr-matrix-banddiag}
For the proof of Theorem \ref{thm:matrix-banddiag}, here we start with the following lemma which is based on the translation of Lemma \ref{lemma:expansion-recursion} by the `matching' used in (\ref{eqn:basis_CSI}): 
\begin{lemma}
\Label{lemma:xjDm-banddiag}
Let $k_0, j, m \in \mathbb{Z}^+$, $\kappa \in\mathbb{Z}$ and $\, \ell _1 := 2m+k_0-\kappa $. 
Under the choices 
{\rm (\ref{eqn:spaces_CSI1})}, 
{\rm (\ref{eqn:spaces_CSI2})},
and {\rm (\ref{eqn:basis_CSI})}, for $\, \kappa \le k_0+m-j $, the function $x^j ({\textstyle \frac{d}{dx}})^m e_n(x)$ can be expressed as a linear combination of $e_{n'}^\Diamond $ $(n'=0,1, ... n+\ell _1)$ 
at most for $n<\ell _1 $, and it can be expressed as a linear combination of $e_{n'}^\Diamond $ 
$(n'=n-\ell _1, \, n-\ell _1 +2, \, n-\ell _1 +4, \, n-\ell _1 +6, \, ...., n+\ell _1)$ for $n\ge \ell _1 $. 
In these linear combinations, all the coefficients are polynomials of ${\ddot{n}}_{k_0,n}$ and $k_0$ with degree not greater than $m$. 
In particular, in the linear combination for $n\ge \ell _1 $, 
with ${\ddot{n}}_{k_0,n}$ defined in {\rm (\ref{eqn:basis_CSI})}, 
the coefficient of the first term with $e_{n-\ell _1 }^\Diamond $ is $\, \displaystyle \left(\frac{i}{2}\right)^{k_0-\kappa -j+m} \left(\frac{1}{2}\right)^j (-1)^m \prod_{t=1}^m \left({\ddot{n}}_{k_0,n}+k_0+t\right)$ when $n+k_0$ is even, and it is $\, \displaystyle \left(-\frac{i}{2}\right)^{k_0-\kappa -j+m} \left(\frac{1}{2}\right)^j \prod_{t=1}^m \left({\ddot{n}}_{k_0,n}-t+1\right)$ when $n+k_0$ is odd.
\end{lemma} \par 
\par The proof is derived directly from Lemma \ref{lemma:expansion-recursion} together with (\ref{eqn:basis_CSI}).

\noindent\noindent{\em Proof of Theorem }\ref{thm:matrix-banddiag}: \quad
\rm From the definition of $s_0$, the inequality ${k_0^\Diamond }\le k_0-s_0$ implies that $ ^\forall m \in \{0,1,2,...,M\}, \,\,\, {k_0^\Diamond }\le k_0+m-\deg p_m$. Hence for every term in the expansion $\displaystyle P(x,{\textstyle \frac{d}{dx}}) = \sum_{m=0}^M \sum_{j=0}^{\deg p_m} p_{m.j}\, x^j ({\textstyle \frac{d}{dx}})^m $, ${k_0^\Diamond }\le k_0+m-j$ holds. Therefore, we can apply Lemma \ref{lemma:xjDm-banddiag} term-wise in this expansion, where $\langle x^j ({\textstyle \frac{d}{dx}})^m e_r, e_n^\Diamond \rangle _{{{\cal H}^\Diamond }} =0$ for 
\par\noindent 
$|r-n|>2m+k_0-{k_0^\Diamond }$ and of course for $|r-n|>2M+k_0-{k_0^\Diamond }$. Hence \par\noindent $b_n^r =\langle Be_r, \, e_n^\Diamond \rangle _{{{\cal H}^\Diamond }} =0$ for $|r-n|>2M+k_0-{k_0^\Diamond }$\, 1.e. (a) holds. 

Next, we will show (b). Since ${\ddot{n}}_{k_0,r} := \lfloor -{\textstyle \frac{k_0+1}{2}} \rfloor + (-1)^{r+k_0+1}\, \lfloor {\textstyle \frac{r+1}{2}} \rfloor$, \par\noindent 
$\displaystyle |{\ddot{n}}_{k_0,r}|\le \frac{r+k_0+4}{2}$. Hence, for fixed $k_0$, for any polynomial $B(x)$, there exists a polynomial $A(x)$ of the same degree as $B(x)$ such that $|B({\ddot{n}}_{k_0,r})|\le A_m(r)$ for $r\in \mathbb{Z}^+$. Since Lemma \ref{lemma:xjDm-banddiag} implies that there exists a polynomial $B_{(n)}(x)$ of degree not greater than $M$ such that $b_n^r=B_{(n)}(r)$ for every $n\in\mathbb{Z}^+$, this fact results in the existence of a polynomial $A(x)$ of degree not greater than $M$ such that $|b_n^r|\le A(r)$ for $r\in \mathbb{Z}^+$,\, i.e. (b) holds. 

Moreover, for $m<M-1$, $\displaystyle\left\langle x^j {\textstyle(\frac{d}{dx})}^m e_r, \, e_{r-(2M+k_0-{k_0^\Diamond })}^\Diamond \right\rangle _{{{\cal H}^\Diamond }} =0$ because 
\par\noindent
 $|r-(r-(2M+k_0-{k_0^\Diamond }))|>2m+k_0-{k_0^\Diamond }$. Hence
\par\noindent 
$\displaystyle b_{r+2M+k_0-{k_0^\Diamond }}^r = \langle Be_r, e_{r+2M+k_0-{k_0^\Diamond }}^\Diamond \rangle _{{{\cal H}^\Diamond }} = \left\langle p_M(x) {\textstyle(\frac{d}{dx})}^M e_r , e_{r-(2M+k_0-{k_0^\Diamond })}^\Diamond \right\rangle _{{{\cal H}^\Diamond }} $. On the other hand, with ${\ddot{n}}_{k_0,r} := \lfloor -{\textstyle \frac{k_0+1}{2}} \rfloor + (-1)^{r+k_0+1}\, \lfloor {\textstyle \frac{r+1}{2}} \rfloor$, Lemma \ref{lemma:xjDm-banddiag} implies that
\begin{eqnarray*}
&& \left\langle x^j {\textstyle(\frac{d}{dx})}^M e_r, \, e_{r-(2M+k_0-{k_0^\Diamond })}^\Diamond \right\rangle _{{{\cal H}^\Diamond }} \\ \\ \displaystyle \,\,\,\,\,\,\,\,\,\,\,\,\,\,\,\, 
&&= 
\left\{\begin{array}{@{\,}ll}\displaystyle (-i)^j \,\left(\frac{i}{2}\right)^{k_0-{k_0^\Diamond }+M} (-1)^M \prod_{t=1}^M \left({\ddot{n}}_{k_0,r}+k_0+t\right)
& \,\,\,\,\, ({\rm if} \,\,\, k_0+r:{\rm even}) \\ \\ \displaystyle 
i^j \,\left(-\frac{i}{2}\right)^{k_0-{k_0^\Diamond }+M} \prod_{t=1}^M \left({\ddot{n}}_{k_0,r}-t+1\right) 
& \,\,\,\,\, ({\rm if} \,\,\, k_0+r:{\rm odd}) 
\end{array}\right. 
\end{eqnarray*}
holds for $r\ge 2M+k_0-{k_0^\Diamond }$, 
because $\displaystyle \left(\frac{1}{2}\right)^j = (\mp i)^j \left(\frac{\pm i}{2}\right)^j$. 
These facts and the relation $\displaystyle \sum_{j=0}^{\deg p_M} p_{M,j} \, (\pm i)^j = p_M(\pm j)$ for $r\ge 2M+k_0-{k_0^\Diamond }$ imply that 
\begin{eqnarray*}
&&\left\langle Be_r, \, e_{r-(2M+k_0-{k_0^\Diamond })}^\Diamond \right\rangle _{{{\cal H}^\Diamond }} \\
&&= 
\left\{\begin{array}{@{\,}ll}\displaystyle p_M(-i) \,\,\left(\frac{i}{2}\right)^{k_0-{k_0^\Diamond }+M} (-1)^M \prod_{t=1}^M \left({\ddot{n}}_{k_0,r}+k_0+t\right)
& \,\,\,\,\, ({\rm if} \,\,\, k_0+r:{\rm even}) \\ \\ \displaystyle 
p_M(i) \,\,\left(-\frac{i}{2}\right)^{k_0-{k_0^\Diamond }+M} \prod_{t=1}^M \left({\ddot{n}}_{k_0,r}-t+1\right) 
& \,\,\,\,\, ({\rm if} \,\,\, k_0+r:{\rm odd}) .
\end{array}\right.
\end{eqnarray*}
>From the definition of ${\ddot{n}}_{k_0,r}$, at least with $r\ge k_0+2M$,\, ${\ddot{n}}_{k_0,r}+k_0+t \le -2$ for \par\noindent $t=0,1,2,...,M$ when $k_0+r$ is even and ${\ddot{n}}_{k_0,r}-t+1 \ge 1$ for $t=0,1,2,...,M$ when $k_0+r$ is odd (where $r=k_0+2M$ is impossible). Since $p_M(\pm i)\ne 0$ from the condition, we have the conclusion $\left\langle Be_r, \, e_{r-(2M+k_0-{k_0^\Diamond })}^\Diamond \right\rangle _{{{\cal H}^\Diamond }} \ne 0$ at least for \par\noindent $r\ge k_0+2M+\max (-{k_0^\Diamond },\, 0)$ i.e. (c) holds.
\hfill\endproof
\par

\end{document}